\documentclass[reqno]{amsart}
\usepackage[utf8]{inputenc}
\usepackage[USenglish]{babel}
\usepackage{amsmath,amsfonts,amssymb,amsthm,dsfont,bm}
\usepackage{enumitem}
\usepackage[usenames,dvipsnames]{xcolor}
\usepackage[colorlinks=true]{hyperref} 
\usepackage{verbatim}
\usepackage{graphicx}
\usepackage{mathabx}
\usepackage{cite}
\hypersetup{linkcolor=Violet,citecolor=PineGreen}
\newtheorem{thm}{Theorem}[section]
\newtheorem{lem}[thm]{Lemma}
\newtheorem{prop}[thm]{Proposition}
\newtheorem{cor}[thm]{Corollary}
\theoremstyle{definition}
\newtheorem{defn}[thm]{Definition}
\newtheorem{rmk}[thm]{Remark}
\newtheorem{rmks}[thm]{Remarks}
\numberwithin{equation}{section}
\definecolor{col}{rgb}{0,0,0.6}
\newcommand{\N}{\mathbb{N}}
\newcommand{\R}{\mathbb{R}}
\newcommand{\Q}{\mathbb{Q}}
\newcommand{\K}{\mathcal{K}}

\newcommand{\T}{\mathcal{T}}

\newcommand{\B}{\mathcal{B}}
\newcommand{\linear}{\mathcal{L}}
\newcommand{\oU}{\overline U}
\newcommand{\ep}{\varepsilon}
\newcommand{\des}{\displaystyle}
\newcommand{\C}{\mathfrak{C}}
\newcommand{\LC}{\mathfrak{LC}}
\newcommand{\SC}{\mathfrak{SC}}

\newcommand{\nbd}{\nobreakdash}

\newcommand{\nti}{{n\to\infty}}

\DeclareMathOperator{\supp}{supp}

\makeatletter
\@namedef{subjclassname@2020}{%
  \textup{2020} Mathematics Subject Classification}

\makeatother
\begin{document}
\title[Topologies of continuity for  Carath\'{e}odory parabolic PDEs]
{Topologies of continuity for  Carath\'{e}odory parabolic PDEs from a dynamical perspective}
\author[I.P.~Longo]{Iacopo P. Longo}
\author[R.~Obaya]{Rafael Obaya}
\author[A.M.~Sanz]{Ana M. Sanz}
\address[Iacopo P. Longo]{Technische Universit\"at M\"unchen,
Forschungseinheit Dynamics,
Zentrum Mathematik, M8,
Boltzmannstra{\ss}e 3,
85748 Garching bei M\"unchen, Germany.}
\email{longoi@ma.tum.de}
\address[R. Obaya]{Departamento de Matem\'{a}tica
Aplicada, Escuela de Ingenier\'{\i}as Indus\-tria\-les (Sede Doctor Mergelina), Universidad de Valladolid,
47011 Valladolid, Spain, and member of IMUVA, Instituto de Investigaci\'{o}n en
Matem\'{a}ticas, Universidad de Va\-lla\-dolid, Spain.}
 \email{rafael.obaya@uva.es}
\address[A.M. Sanz]{Departamento de Did\'{a}ctica de las Ciencias Experimentales, Sociales y de la Matem\'{a}tica,
Facultad de Educaci\'{o}n, Universidad de Valladolid, 34004 Palencia, Spain,
and member of IMUVA, Instituto de Investigaci\'{o}n en  Mate\-m\'{a}\-ti\-cas, Universidad de
Valladolid.} \email{anamaria.sanz@uva.es}
\thanks{All authors were partly supported by MICIIN/FEDER project
RTI2018-096523-B-I00 and by the University of Valladolid under project PIP-TCESC-2020.  I.P.~Longo was also partly supported by the European Union’s Horizon 2020 research and innovation programme under the Marie Skłodowska-Curie grant agreement No 754462 and by TUM International Graduate
School of Science and Engineering (IGSSE)}
\date{}
\begin{abstract}
Systems of non-autonomous parabolic partial differential equations over a bounded domain with nonlinear term of Carath\'{e}odory type are considered. Appropriate topologies on sets of Lipschitz Carath\'{e}odory maps are defined in order to have a continuous dependence of the mild solutions with respect to the variation of both the nonlinear term and the initial conditions, under different assumptions on the bound-maps of the nonlinearities.
\end{abstract}
\maketitle
\section{Introduction}\label{secintro}
The study of the topologies of continuity for Carath\'{e}odory ordinary differential equations (ODEs for short) is a classical question with important implications in the field of non-autonomous differential equations and  in dynamical systems, and their applications in mathematical modelling. In particular, Artstein~\cite{paper:ZA1,paper:ZA2,paper:ZA3}, Heunis~\cite{paper:AJH}, Miller and Sell~\cite{book:RMGS,paper:RMGS1}, Neustadt~\cite{paper:LWN}, Sell~\cite{paper:GS1,book:GS}, and some more references therein, introduced and studied strong and weak topologies of integral type on spaces of Lipschitz Carath\'{e}odory functions. These Carath\'{e}odory vector fields define a set $E$, and the convergence of a sequence $(f_n)_{n \geq 1}$ with respect to these topologies in $E$ requires the convergence of the integral of the evaluation of the functions $f_n$, either pointwise in $\R^N$ (topologies $\T_D$ and $\sigma_D$), or uniformly on any bounded set of continuous functions (topology $\T_B$). Recently, Longo et al.~\cite{paper:LNO1,paper:LNO2}  completed some parts of this theory by introducing the strong and weak topologies $\T_\Theta$ and $\sigma_\Theta$, where $\Theta$ is a suitable set of moduli of continuity. The motivation is that when the set of $m$-bounds of the functions in $E$ is equicontinuous, then given $j \in \N$ and a compact interval $I \subset \R$, the  solutions of the equations for $f\in E$ defined on $I$ and bounded by $j$  admit a common modulus of continuity and are thus included in a compact subset $\K^I_j$ of continuous functions. Therefore, it is possible to define a countable set of moduli of continuity $\Theta$ so that the topologies $\T_\Theta$ and $\sigma_\Theta$, defined by the uniform convergence of the integral of the evaluation of the functions of $E$ on each compact set $\K^I_j$ are of continuity, i.e., if the vector fields $(f_n)_{n\geq 1}$ converge to $f$ in E, and the initial data $(x_n)_{n\geq 1}$ converge to $x$ in $\R^N$, then the solutions of the differential equations $y(t,f_n,x_n)$ converge to $y(t,f,x)$ uniformly on compact time intervals.

The references Longo et al.~\cite{paper:LNO3,paper:LNO4} provide an extension of these methods to  Carath\'{e}odory delay differential equations. To overcome the problem that  in general the solutions of the equations do not share the same modulus of continuity as the initial data,  it is necessary to introduce hybrid topologies from suitable pairs of topologies of the types defined above which treat the present and the past components of the equations differently.  In addition,  the conclusions on the new topologies of continuity for  delay differential equations are applied to develop dynamical methods to investigate non-autonomous models in populations dynamics.

In this paper we introduce  appropriate versions of these topologies on sets $E$ of measurable functions $f:\R\times \oU \times \R^N \rightarrow \R^N$ with a Lipschitz variation on the state component, defining Carath\'{e}odory parabolic partial differential equations (PDEs for short) on an open and bounded domain $U\subset \R^M$ with a smooth boundary,   and prove the continuous variation of the mild solutions of the associated abstract Cauchy problems (ACPs for short) in the space $C(\oU,\R^N)$ with respect to the initial condition and the function  $f \in E$. Although this theory follows the same ideas explained above for the case of ordinary and delayed  differential equations, new  nontrivial technical challenges appear in this setting, as a consequence of the unbounded behavior of the realization of the Laplacian subject to the different boundary conditions.

The paper is organized in four sections. In Section~\ref{secCaratheodoryPDEs} we consider non-autonomous $N$-dimensional parabolic initial boundary value (IBV for short) problems
\begin{equation*}
\left\{\begin{array}{l} \des\frac{\partial y}{\partial t}  =
 \Delta \, y+f(t,x,y)\,,\quad t>0\,,\;\,x\in U,
  \\[.2cm]
By=0\,,\quad  t>0\,,\;\,x\in \partial U,
\\
y(0,x)=z(x)\,,\quad  x\in \oU,
\end{array}\right.
\end{equation*}
with Neumann, Robin or Dirichlet boundary conditions, where the nonlinear term $f(t,x,y)$ belongs to the space ${\LC}$ of Lipschitz Carath\'{e}odory functions. We denote by $B_j$ the closed ball of $\R^N$ centered at the origin and by ${\B}_j=C(\oU,B_j)$. We consider the  infinite dimensional Banach space $X=C({\oU},\R^N)$ if the boundary conditions are of Neumann or Robin type whereas  we will take the space $X=C_0({\oU},\R^N)$ if the boundary conditions are of Dirichlet type. We  transform the former problem into an ACP in  the space $X$ and prove, for each $z \in X$, the existence  of a unique mild solution $u(t,f,z)$ satisfying $u(0,f,z)=z$.

In Section~\ref{sectopo} the appropriate versions of the classical topologies ${\T_B}, {\T_D}, \sigma_{D}$ on the space ${\SC}$ of strong Carath\'{e}odory maps, as well as the new topologies in this PDEs context $\T_{\widetilde DD}$ and $\sigma_{\widetilde DD}$,  are introduced. In the definition of the topologies ${\T}_\Theta$ and $\sigma_\Theta$, a suitable set of moduli of continuity $\Theta=(\Theta_1,\Theta_2)$ is now formed by a pair of parametric families of ordered moduli of continuity that, for each compact interval $I \subset \R$ and each $j \in \N$,  determine bounded sets $\mathcal{H}^I_j \subset C(I,{\B}_j)$  with a precise equicontinuous variation on $I$, and the seminorms are  given by the supremum of the integrals of the evaluation of the functions in ${\SC}$ along the maps in $\mathcal{H}^I_j$. The reason for this is that the continuity properties that the mild solutions inherit from the equations involve the regularity properties of the semigroup of operators $(e^{tA})_{t\geq 0}$ generated by the realization of the Laplacian in $X$. Then, for a set $E\subset \LC$ we study the equivalence of some of the previous topologies under different assumptions on the $m$-bounds and/or the $l$-bounds of the maps in $E$.


Finally, Section~\ref{secContMildSolutions} states precise conditions implying that  the previous topologies are topologies of continuity for the parabolic PDEs. Associated to a family $E \subset {\LC}$ with $L^1_{loc}$-equicontinuous $m$-bounds we determine a suitable set of moduli of continuity $\Theta$ and a suitable set of radii $\mathcal{R}$ that provide sets ${\K}^I_j  \subset \mathcal{H}^I_j$   satisfying that given a compact interval $I \subset \R$ and $j \in \N$,  if $u(t,f,z)$ is a family of mild solutions of the Carath\'{e}odory ACPs defined on $I$ for each $f \in E$ and $\|u(t,f,z)\|\leq j$ for all $t\in I$, then $u(\cdot,f,z)\in \K_j^I$. Thus, the topologies $\T_{\Theta\mathcal{R}}$ and $\sigma_{\Theta\mathcal{R}}$ whose seminorms are  given by the supremum of the integral of the evaluation of the functions on the  sets ${\K}^I_j$ are both of them topologies of continuity for the parabolic equations. Since $\T_{\Theta\mathcal{R}} \leq {\T}_{\Theta}$ and  $\sigma_{\Theta\mathcal{R}}\leq \sigma_{\Theta}$, it is obvious that ${\T}_{\Theta}$ and  $\sigma_{\Theta}$ is a pair of topologies of continuity as well. On the other hand, assuming that the family $E$ admits  $L^1_{loc}$-bounded $l$-bounds, it is proved the continuous variation of the mild solutions $u(t,f,z)$ when $f$ varies on E with  respect to any strong topology $\T_D$, $\T_\Theta$ or $\T_B$, all of them equivalent in this case. Finally, the weaker topologies of pointwise type $\T_{\widetilde DD}$ and $\sigma_{\widetilde DD}$ are also topologies of continuity in the previous situations, if in addition $E$  has $L^1_{loc}$-bounded $l$-bounds with respect to a modulus of continuity $\widehat \theta$ in the variable $x$.
\section{Mild solutions for Carath\'{e}odory parabolic PDEs}\label{secCaratheodoryPDEs}
In this section we specify the kind of $N$-dimensional parabolic IBV
problems for $y(t,x)$, with Neumann, Robin or Dirichlet boundary conditions, under consideration:
\begin{equation}\label{eq:pde}
\left\{\begin{array}{l} \des\frac{\partial y_i}{\partial t}  =
 \Delta \, y_i+f_i(t,x,y)\,,\quad t>0\,,\;\,x\in U,\; 1\leq i\leq N,
  \\[.2cm]
By_i:=\alpha_i(x)\,y_i+\kappa\,\des\frac{\partial y_i}{\partial n} =0\,,\quad  t>0\,,\;\,x\in \partial U,\; 1\leq i\leq N,
\\[.2cm]
y_i(0,x)=z_i(x)\,,\quad  x\in \oU,\; 1\leq i\leq N,
\end{array}\right.
\end{equation}
where  $U$, the spatial domain, is a bounded, open and
connected  subset of $\R^M$ ($M\geq 1$) with a sufficiently smooth boundary
$\partial U$; $\Delta$ is the Laplacian operator on $\R^M$ and the nonlinear terms are given by a map $f:\R\times \oU\times \R^N\to \R^N$ with components $(f_1,\ldots,f_N)$ within a class of so-called {\em Carath\'{e}odory maps\/}, namely, the class of Lipschitz Carath\'{e}odory maps. The problem has Dirichlet boundary conditions if $\kappa=0$ and $\alpha_i(x)\equiv 1$ for $1\leq i\leq N$; Neumann boundary conditions if $\kappa=1$ and  $\alpha_i(x) \equiv 0$ for $1\leq i\leq N$; and Robin boundary conditions if  $\kappa=1$ and $\alpha_i:\partial U\to \R$ is a nonnegative  sufficiently regular map for $1\leq i\leq N$. Recall that $\partial/\partial n$ denotes the outward normal derivative at the  boundary. Finally, the initial condition is given by the values of a map $z=(z_1,\ldots,z_N):\oU\to \R^N$, which we choose to be continuous, with the compatibility condition $z(x)=0$ for all $x\in \partial U$ in the Dirichlet case.

We define the standard classes of Carath\'{e}odory maps. Let $L^1_{loc}$ denote the space of  locally integrable real maps defined on $\R$ and let   $|\cdot|$ denote the norm on the $N$\nbd-dimensional space  $\R^N$. For convenience, we take $|y|=|y_1|+\ldots +|y_N|$.  The symbol  $B_r$ will denote the closed ball of $\R^N$ centered at the origin with radius $r$.

We will consider, and denote by $\C$, the set of Carath\'{e}odory functions $f\colon\R\times\oU\times\R^N\to \R^N$, $(t,x,y)\mapsto f(t,x,y)$ satisfying
\begin{enumerate}[label=\upshape(\textbf{C\arabic*}),leftmargin=27pt,itemsep=2pt]
\item\label{C1}  $f$ is Borel measurable and
\item\label{C2} for every compact set $K\subset\R^N$ there exists a nonnegative function $m^K\in L^1_{loc}$, called \emph{$m$-bound}, such that for almost every $t\in\R$ one has $|f(t,x,y)|\le m^K(t)$ for all $x\in\oU,\, y\in K$.
\end{enumerate}
\smallskip
We introduce the classes  of Carath\'eodory functions which are subsequently used. In all the cases, we identify the maps within a certain class which differ on a negligible set.
\begin{defn}\label{def:SC}
A function $f\colon\R\times\oU\times\R^N\to \R^N$, $(t,x,y)\mapsto f(t,x,y)$  is said to be \emph{strong Carath\'eodory}, and we will write $f\in \SC$, if it satisfies~\ref{C1},~\ref{C2} and
\begin{enumerate}[label=\upshape(\textbf{S}),leftmargin=27pt,itemsep=2pt]
\item\label{S} for almost every $t\in\R$, the function $f(t,\cdot,\cdot)$ is continuous.
\end{enumerate}
\end{defn}
\begin{rmk}\label{rmk:hypotheses}
Note that, if $f$ is Borel measurable and it satisfies \ref{S}, then (i) $f$ is Lebesgue measurable in $t$ for each fixed $(x,y)\in \oU\times \R^N$ (see, e.g., Lemma~5.1.2 in Cohn~\cite{book:Cohn}) and (ii) it is continuous in $(x,y)$ for almost every fixed $t\in\R$. In fact, the hypotheses are often formulated by requiring (i) and (ii). But using standard arguments of measure theory one can check that, if (i) and (ii) hold, then there exists a Borel measurable map $g\colon\R\times\oU\times\R^N\to \R^N$ satisfying~\ref{S}, and there exists a set $N\subset\R$ of null measure such that $f$ coincides with $g$ on the set $(\R\setminus N)\times \oU\times\R^N$. To give the idea, the map $g$ is built as a pointwise limit of a sequence of Borel measurable maps  $f_n:=f\,$\raisebox{2pt}{$\chi$}$_{G_n\times \oU \times B_n}$ for an expanding sequence of Borel sets $G_n$ of $\R$ such that $\R\setminus (\cup_{n\geq 1} G_n)$ has null measure and the restriction of $f$ to $G_n\times \oU \times B_n$ is continuous for every $n\geq 1$.
\end{rmk}
\begin{defn}\label{def:LC}
A function $f\colon\R\times\oU\times\R^N\to \R^N$, $(t,x,y)\mapsto f(t,x,y)$  is said to be \emph{Lipschitz Carath\'eodory}, and we will write $f\in \LC$, if $f\in\SC$ and
\begin{enumerate}[label=\upshape(\textbf{L}),leftmargin=27pt,itemsep=2pt]
\item\label{L} for every compact set $K\subset\R^N$ there exists a nonnegative function $l^K\in L^1_{loc}$,  called \emph{$l$-bound}, such that for almost every $t\in\R$, $|f(t,x,y_1)-f(t,x,y_2)|\le l^K(t)\,|y_1-y_2|$ for all $x\in\oU,\,y_1,y_2\in K$.
\end{enumerate}
\end{defn}
In particular, for each compact set $K\subset\R^N$, we refer to \emph{the optimal $m$-bound} and \emph{the optimal $l$-bound} of $f$ as to
the maps defined almost everywhere by
\begin{equation}\label{eqOptimalMLbound}
\begin{split}
m^K(t)&=\sup_{x\in \oU,\,y\in K}|f(t,x,y)|\quad \mathrm{and}
\\
\quad l^K(t)&=\sup_{\substack{x\in \oU,\,y_1,y_2\in K\\[3pt] y_1\neq y_2}}\frac{|f(t,x,y_1)-f(t,x, y_2)|}{|y_1-y_2|}\, ,
\end{split}
\end{equation}
respectively. Clearly, for each compact set $K\subset\R^N$ the suprema in~\eqref{eqOptimalMLbound}  can be taken for some countable dense subsets of $K$ and $\oU$, respectively, leading to the same definition, which guarantees that the functions defined in~\eqref{eqOptimalMLbound} are measurable. Moreover, we consider a subclass within the class of Lipschitz Carath\'eodory maps, by fixing a precise continuous behaviour with respect to $x$, e.g., H\"{o}lder continuity or a Lipschitz character.
\begin{defn}\label{def:LC con modulo en x}
Given a modulus of continuity, i.e.,~a  non-decreasing continuous function  $\widehat \theta:\R^+ \to\R^+$ with $\widehat \theta(0)=0$, a function $f\colon\R\times\oU\times\R^N\to \R^N$, $(t,x,y)\mapsto f(t,x,y)$  is said to be \emph{Lipschitz Carath\'eodory with respect to the modulus of continuity $\widehat\theta$ in the variable $x$}
 if it satisfies~\ref{C1},~\ref{C2}, and
\begin{enumerate}[label=\upshape(\textbf{L}$_{\smash{\widehat\theta}}$), leftmargin=27pt,itemsep=2pt]
\item\label{LT} for every compact set $K\subset\R^N$ there is a nonnegative function $l^K\in L^1_{loc}$ such that for almost every $t\in\R$,
\begin{equation*}
|f(t,x_1,y_1)-f(t,x_2,y_2)|\le l^K(t)\big(\widehat\theta(|x_1-x_2|)+|y_1-y_2|\big)
\end{equation*}
for all $x_1,x_2\in\oU,\,y_1,y_2\in K$.
\end{enumerate}
\end{defn}
Note that these maps also satisfy~\ref{S}. In particular, for each compact set $K\subset\R^N$, we refer to \emph{the optimal $l$-bound of $f$ with respect to $\widehat \theta$} as to the map defined almost everywhere by
\begin{equation}
 l^K(t)=\sup_{\substack{x_1,x_2\in \oU,\, y_1,y_2\in K\\[3pt] x_1\neq x_2,\,y_1\neq y_2}}\frac{|f(t,x_1,y_1)-f(t,x_2, y_2)|}{\widehat\theta(|x_1-x_2|)+|y_1-y_2|}\, .
\label{eq:Lbound-theta}
\end{equation}

The aim of this section is to study the existence and uniqueness of mild solutions for the Carath\'{e}odory ACPs in an appropriate Banach space associated to the parabolic problems~\eqref{eq:pde}.
\begin{rmk}\label{rmk:notation X}
Trying to keep a common notation, we will consider the space $X=C(\oU,\R^N)$ of the continuous functions on $\oU$ taking values in $\R^N$ if the boundary conditions are of Neumann or Robin type, whereas we will take the space $X=C_0(\oU,\R^N)$ of the continuous functions on $\oU$ vanishing on the boundary $\partial U$  if the boundary conditions are of Dirichlet type.
\end{rmk}

Using for simplicity the same notation, when dealing with the laplacian with Neumann or Robin boundary conditions,  for each component $i=1,\ldots,N$ we consider on the Banach space $Y=C(\oU)$ of real continuous maps on $\oU$ endowed with the sup-norm, the differential operator $A_i^0 z_i = \Delta z_i$ with domain $D(A_i^0)$ given by
\[
\left\{z_i\in C^2(U)\cap C^1(\oU)\;\Big|\; A_i^0z_i\in C(\oU)\,,\;
\alpha_i(x)\,z_i(x)+\des\frac{\partial z_i}{\partial n}(x)=0\, \;\forall \,x\in \partial U \right\},
\]
whereas if the boundary conditions are of Dirichlet type, then we take the Banach space $Y=C_0(\oU)$  and the domain  $D(A_i^0)=\big\{z_i\in C^2(U)\cap C_0(\oU)\;|\; A_i^0z_i\in C_0(\oU)\big\}$.

In both cases  the closure $A_i$ of $A_i^0$ in  $Y$ is a sectorial operator which generates an
analytic semigroup of bounded  linear operators $(e^{tA_i})_{t\geq 0}\subset \mathcal{L}(Y)$, which is strongly continuous (that is, $A_i$ is densely defined),
and $e^{tA_i}$ is compact for all $t>0$ (see, e.g., Smith~\cite{book:smit}). Furthermore, using the theory of positive  semigroups of operators, it is deduced that $\|e^{tA_i}\|\leq 1$ for all $t\geq 0$ (see~\cite[Corollary 7.2.4]{book:smit}).

At this point, we can consider $A=\Pi_{i=1}^N  A_i$ which is a sectorial operator with domain $D(A)=\Pi_{i=1}^N D(A_i)$ on the product Banach space  $X=Y^N$ endowed with the norm $\|(z_1,\ldots,z_n)\|_1=\sum_{i=1}^N\|z_i\|$. Then, $\{e^{tA}=\Pi_{i=1}^N e^{tA_i}\mid t\geq 0\}$ defines a positive semigroup of operators which satisfies $\|e^{tA}\|\leq 1$ for all $t\geq 0$  and $e^{tA}$ is compact for all $t>0$.

For convenience, hereafter we consider the sup-norm  on the space $X$, denoted  by $\|z\|=\sup_{x\in \oU} |z(x)|$. As a consequence of the bound in the previous paragraph, now for the operator norm on $\mathcal{L}(X)$ associated to the sup-norm, we have that
\begin{equation}\label{eq:bound semigroup}
\|e^{tA}\|\leq N\quad \text{for all}\; t\geq 0\,.
\end{equation}
We will be using this bound later.

Given  $f\in\SC$, for each $t\in\R$ and $z\in X$, let  $\tilde f(t,z):\oU\to\R^N$ be defined by
\begin{equation}\label{eq:tilde f}
\tilde f(t,z)(x)=f(t,x,z(x))\,,  \quad x\in \oU.
\end{equation}
With Dirichlet boundary conditions we will assume with no further mention that
\begin{equation*}
f(t,x,0)=0  \quad \text{for all}\;t\in \R\,,\;x\in \partial U,
\end{equation*}
so that for all $t\in\R$ and $z\in X=C_0(\oU,\R^N)$, $\tilde f(t,z)(x)=0$ for all $x\in\partial U$, and provided that $\tilde f(t,z)$ is continuous on $\oU$, then $\tilde f(t,z)\in C_0(\oU,\R^N)$. Note that, since $f$ satisfies~\ref{S}, for almost every $t\in\R$ and for all $z\in X$, the map $\tilde f(t,z)\in X$. We now list some properties of the map $\tilde f$.
\begin{prop}\label{prop:f tilde medible}
  Let $f\in\SC$ and consider the map $\tilde f:\R\times X \to X$ defined for almost every $t\in\R$ and every $z\in X$ by~\eqref{eq:tilde f}. Then, the following properties hold:
  \begin{itemize}
    \item[\rm{(i)}] There exists a Borel measurable map $\tilde g:\R\times X\to X$ and a set of null measure $N\subset \R$ such that  $\tilde f = \tilde g$ on the set $(\R\setminus N)\times X$.
    \item[\rm{(ii)}] If $m^j\in L^1_{loc}$ is an $m$-bound of $f$ on the closed ball $B_j$, then for almost every $t\in\R$,
\begin{equation}\label{eq:bound}
\|\tilde f(t,z)\|\leq m^j(t)\in L^1_{loc}\quad \text{for all}\; z\in X\,,\; \|z\|\leq j\,.
\end{equation}
    \item[\rm{(iii)}] If $I\subseteq \R$ is an interval and $u\in C(I,X)$, then the map defined almost everywhere $I\to X$, $t\mapsto\tilde f(t,u(t))$ is measurable on $I$.
        \item[\rm{(iv)}] If in addition $f\in \LC$, and $l^j\in L^1_{loc}$ is an $l$-bound of $f$ on $B_j$, then for almost every $t\in\R$, $\|\tilde f(t,z_1)- \tilde f(t,z_2)\|\leq l^j(t)\,\|z_1-z_2\|$ for all $z_1, z_2\in X$ with $\|z_1\|, \|z_2\|\leq j$.
  \end{itemize}
\end{prop}
\begin{proof}
As noted in Remark~\ref{rmk:hypotheses},  one can consider the Borel measurable map $g$ obtained as the pointwise limit of a sequence of Borel measurable maps $f_n=f\,$\raisebox{2pt}{$\chi$}$_{G_n\times \oU \times B_n}$ for an expanding sequence of sets $G_n$ of $\R$ such that $\R\setminus (\cup_{n\geq 1} G_n)$ has null measure and the restriction of $f$ to $G_n\times \oU \times B_n$ is continuous for every $n\geq 1$. It is not difficult to check that if $\B_n:=\{z\in X\mid \|z\|\leq n\}$, then  $\tilde f$ restricted to the Borel measurable set $G_n\times \B_n$ is also continuous. Then,   $\tilde f_n:= \tilde f \,$\raisebox{2pt}{$\chi$}$_{G_n\times \B_n}$, $n\geq 1$ are  Borel measurable maps and their pointwise limit $\tilde g$ is a Borel measurable map which satisfies $\tilde f = \tilde g$ on the set $(\cup_{n\geq 1} G_n)\times X$, so that (i) holds.

For (ii), note that by~\ref{C2},  if $\|z\|\leq j$, then for almost every $t\in\R$,  $|f(t,x,z(x))|\leq m^j(t)$ for all $x\in\oU$. From here,~\eqref{eq:bound} holds. Finally,~(iii) follows from (i) and the continuity of $u$ on $I$, and~(iv) follows from~\ref{L} straightaway.
\end{proof}
We are in a position to consider a  Carath\'{e}odory ACP associated to the parabolic problem~\eqref{eq:pde} with $f\in\SC$, that is, an evolution equation in the Banach space $X$ (see Remark~\ref{rmk:notation X})  of Carath\'{e}odory type:
\begin{equation}\label{eq:acp}
\left\{\begin{array}{l} u'(t)  =
 A\, u(t)+\tilde f(t,u(t))\,,\quad t>0\,,\\
u(0)=z\,.
\end{array}\right.
\end{equation}
The fact that $A$ is an unbounded operator on $X$ prevents the application of the basic theory developed by Holly and Orewczyk~\cite{paper:HO}  for Carath\'{e}odory abstract ODEs. In any case, one cannot hope for more than solutions in the sense of Carath\'{e}odory.
It is standard to search for solutions within the class of the so-called {\it mild solutions}, that is, continuous solutions of the associated integral equation. This integral equation is a generalization of the variation of constants formula in standard ODEs' linear Cauchy problems to the present abstract setting. It is important to note that, since the variation of $\tilde f(t,u(t))$ in $t$ is just measurable, we must deal with the Bochner integral (see, e.g., the classical reference Hille and Phillips~\cite{book:HP}), which is an extension of the Lebesgue's integration theory to vector-valued functions. A summary of this theory can be found in Engel and Nagel~\cite[Appendix C]{book:EN}. One basic issue  is that, if the semigroup of operators $(e^{tA})_{t\geq 0}$ is strongly continuous, as it is, and $\tilde f(t,u(t))$ is measurable,  the integrand in~\eqref{eq:variation constants} is measurable, which is a pre-requisite for the integral to make sense.
\begin{defn}\label{def:mildsolution}
A map $u\in C([0,\delta], X)$ is a  {\em mild solution\/} of~\eqref{eq:acp}  on $[0,\delta]$  
if it satisfies
\begin{equation}\label{eq:variation constants}
u(t)=e^{tA}\,u(0) +\int_0^t e^{(t-s)A}\,\tilde f(s,u(s))\,ds\,,\quad t\in[0,\delta]\,.
\end{equation}
\end{defn}

Being not aware of any appropriate reference for this setting, we offer the basic theory on the existence and uniqueness of mild solutions for the Carath\'{e}o\-do\-ry ACP~\eqref{eq:acp} following the approach in Lunardi~\cite{book:luna} in the case of more regular problems, namely, the case in which $\tilde f$ is continuous and satisfies a Lipschitz condition in the second variable. Here we lack the continuity and, in order to have a sufficient Lipschitz-like condition on $\tilde f$  to guarantee uniqueness of solutions, we  restrict ourselves to nonlinear terms $f$ in the class of Lipschitz Carath\'{e}odory maps.
\begin{rmks}\label{rmk:puntoPartida}
1. Note that everything that we are doing for the IBV problem~\eqref{eq:pde} and the associated ACP~\eqref{eq:acp} can be done for a problem
defined on a time interval $[a,\infty)$ for each $a\in \R$, just with some small modifications.  The problems have been written down with initial condition at $t=0$ for the sake of simplicity.

2. It is not difficult to check that if $u\in C([a,T], X)$ is a mild solution of the ACP~\eqref{eq:acp} with initial condition at $a$, then $u$ is also a mild solution on any smaller interval $[t_0,T]$ with $a<t_0<T$, by taking $u(t_0)$ as the initial value at $t_0$. We just mention that one needs to interchange the continuous operator $e^{(t-t_0)A}$, $t_0\leq t\leq T$ with a Bochner integral and this is justified by, e.g.,  Theorem~3.7.12 in~\cite{book:HP}.
\end{rmks}

We include the proof of the following result for the sake of completeness, but the main arguments are taken from~\cite[Theorem~7.1.2]{book:luna}.
\begin{thm}\label{thm:existence mild}
Let $f\in \LC$. Then, given any $z_0\in X$ there exist $r, \delta>0$ such that, if $z\in X$ with $\|z-z_0\|\leq r$, problem~\eqref{eq:acp}  has a unique mild solution $u$ in the space $C([0,\delta],X)$.
\end{thm}
\begin{proof}
The proof of existence is based on the contraction mapping theorem. Take an integer $j\geq 4\,N\,\|z_0\|$ in such a way that, if $\|z-z_0\|\leq j /(4\,N)=:r $, then $\sup_{t\geq 0}\|e^{tA}\,z\|\leq j/2$ (see~\eqref{eq:bound semigroup}). Let $m^j(t)$ (resp.~$l^j(t)$) be the optimal $m$-bound (resp.~$l$-bound) of $f$ on $B_j$ and consider the normed space
\[
Y_j=\{u\in C([0,\delta],X)\mid \|u(t)\|\leq j\;\text{for all}\;t\in[0,\delta]\}\,,
\]
with the sup-norm, where $\delta>0$ is to be determined. Note that $Y_j$ is the closed ball in $C([0,\delta],X)$ of radius $j$. Fixed $z\in X$ with $\|z-z_0\|\leq r$, define the operator $\Gamma:Y_j\to Y_j$ by
\[
\Gamma(u)(t):= e^{tA}\,z +\int_0^t e^{(t-s)A}\,\tilde f(s,u(s))\,ds\,,\quad t\in[0,\delta]\,.
\]
First, since $(e^{tA})_{t\geq 0}$ is strongly continuous, $t\in [0,\delta]\mapsto e^{tA}\,z\in X$ is continuous. Then, using~\eqref{eq:bound semigroup}, \eqref{eq:bound} and~\ref{S}, a combination of the theorem of continuity under the integral sign and the absolute continuity of a Bochner integrable map
 guarantees that $\Gamma(u)\in C([0,\delta],X)$.  To see that $\Gamma(u)\in Y_j$ and that $\Gamma$ is a contractive map, we need to find the appropriate value of $\delta$. Thanks to the absolute continuity of the Lebesgue integral, we can choose a $\delta>0$ such that $\int_0^{\delta} m^j(s)\,ds\leq j/(2\,N)$ and $\int_0^{\delta} l^j(s)\,ds\leq 1/(2\,N)$. Then, for $u\in Y_j$, by Proposition~\ref{prop:f tilde medible}(ii),
\begin{align*}
\|\Gamma(u)\|&\leq  \sup_{t\in [0,\delta]}\|e^{tA}\,z\| + \sup_{t\in [0,\delta]}  \left\|\, \int_0^t e^{(t-s)A}\,\tilde f(s,u(s))\,ds\,\right\|  \\
    &\leq   \frac{j}{2}+ N\int_0^\delta m^j(s)\,ds\leq  \frac{j}{2}+  \frac{j}{2}=j \,
\end{align*}
and, for all $u_1, u_2\in Y_j$, this time by Proposition~\ref{prop:f tilde medible}(iv),
\begin{align*}
\|\Gamma(u_1)-\Gamma(u_2)\| &= \sup_{t\in [0,\delta]} \left\|\, \int_0^t e^{(t-s)A}\,\big(\tilde f(s,u_1(s))- \tilde f(s,u_2(s)) \big)\,ds  \,\right\| \\
    &\leq  N\, \|u_1-u_2\|\int_0^\delta l^j(s)\,ds\leq \frac{1}{2}\,\|u_1-u_2\|\,.
\end{align*}
By the contraction mapping theorem, there exists a unique fixed point of $\Gamma$ on $Y_j$. In particular this means that there exists a  mild solution on $[0,\delta]$ for the Carath\'{e}odory ACP~\eqref{eq:acp}.

To prove uniqueness, if there were two mild solutions $u_1$ and $u_2$ in $C([0,\delta],X)$, define $t_0=\sup\{t\in [0,\delta]\mid u_1(s)=u_2(s) \text{ for } 0\leq s\leq t \}$ and let $z_1=u_1(t_0)=u_2(t_0)$. Assume by contradiction that $t_0<\delta$ and note that by the previous reasoning, the problem
\begin{equation*}
\left\{\begin{array}{l} u'(t)  =
 A\, u(t)+\tilde f(t,u(t))\,,\quad t>t_0\,,\\
u(t_0)=z_1\,
\end{array}\right.
\end{equation*}
has a unique mild solution in a set
\[
\tilde Y=\{u\in C([t_0,t_0+\ep],X)\mid \|u(t)\|\leq k\;\text{for all}\;t\in[t_0,t_0+\ep]\}\,,
\]
for large enough $k$ and small enough $\ep>0$.  Now, both $u_1$ and $u_2$ are bounded on $[t_0,\delta]$ and, as recalled in Remark~\ref{rmk:puntoPartida}.2,  they are mild solutions of the previous problem on $[t_0,t_0+\ep]$ as far as $t_0+\ep<\delta$, so that $u_1,u_2\in\tilde Y$ for an appropriate choice of $k$ and $\ep$. Since, by the definition of $t_0$,  they are different, we find a contradiction.
The proof is finished.
\end{proof}


In what refers to the maximally defined solution, for each $f\in \LC$ and $z\in X$ we define $\beta(f,z)=\sup\{\delta>0\mid \exists \text{ a mild solution } u_\delta \text{ on } [0,\delta] \}$ and $u(t):=u_\delta(t)$ if $t\leq \delta$. Then, by the uniqueness result, $u$ is well-defined on the interval $I_{f,z}:=\bigcup\{[0,\delta]\mid \exists \text{ a mild solution } u_\delta \text{ on } [0,\delta] \}$ and $\beta(f,z)=\sup I_{f,z}$. Then, we have the following result on the behaviour of the mild solution which resembles the well-known result for standard ODEs.
\begin{prop}
Let $f\in \LC$. Assume that  $z\in X$ is such that $I_{f,z}=[0,T)$ with $T<\infty$. Then, $t\mapsto \|u(t)\|$ is unbounded on $I_{f,z}$.
\end{prop}
\begin{proof}
Assume by contradiction that $u$ is bounded on $[0,T)$. Then, since $u(t)$ satisfies~\eqref{eq:variation constants} on $[0,T)$ and for an appropriate $j\geq 1$ and for almost every $s\in [0,T]$, $\|\tilde f(s,u(s))\|\leq m^j(s)\in L^1([0,T])$, $u$ can be extended with continuity to the closed interval $[0,T]$. But then it is not difficult to check that we can extend the mild solution beyond $T$, by gluing $u(t)$ on $[0,T]$ to the mild solution of the ACP
\begin{equation*}
\left\{\begin{array}{l} v'(t)  =
 A\, v(t)+\tilde f(t,v(t))\,,\quad t>T\,,\\
v(T)=u(T)\,
\end{array}\right.
\end{equation*}
which is defined at least on an interval $[T,T+\delta]$, by Theorem~\ref{thm:existence mild}. This is a contradiction with the fact that $I_{f,z}=[0,T)$ was the maximal interval of definition of the mild solution for $f$ and $z$. The proof is finished.
\end{proof}
\section{Strong and weak topologies and topological properties of the bound-maps}\label{sectopo}
We endow the space $\SC$ with suitable strong and weak topologies. The terminology that we use is parallel to the one used in the works by Longo et al.~\cite{paper:LNO1,paper:LNO2}. As a rule, when inducing a topology on a subspace, we will denote the induced topology with the same symbol. Firstly, we consider some integral-like topologies which are the natural adaptation of the classical ones in the context of Carath\'{e}odory ODEs, then we introduce two new topologies $\T_{\widetilde DD}$ and $\sigma_{\widetilde DD}$ in this PDEs context, and finally we define the topologies $\T_{\Theta}$ and $\sigma_{\Theta}$ associated to adequate sets of moduli of continuity.

As explained in Section~\ref{secCaratheodoryPDEs}, linked to the IBV problem \eqref{eq:pde} one builds the Carath\'{e}odory ACP~\eqref{eq:acp} in the Banach space $X$ with the sup-norm $\|\cdot\|$ (see Remark~\ref{rmk:notation X}). For this reason, in the definitions it naturally appears the sup-norm of maps defined on $\oU$ with values in  $\R^N$, where it used to appear the norm in $\R^N$ for standard ODEs. To this end, recall the definition of $\tilde f(t,z)$ given in~\eqref{eq:tilde f}. Also, given  $f\in\SC$, Proposition~\ref{prop:f tilde medible} justifies that all the following definitions make sense.

To simplify the writing, we denote by $\B_j$ the closed ball in $C(\oU,\R^N)$ centered at the null map with radius $j$, that is, $\B_j=C(\oU,B_j)$, and whenever we write  $I=[q_1,q_2],\,  q_1,q_2\in\Q$, it is understood that $q_1<q_2$.
\begin{defn}[Topology $\T_{B}$]
We call $\T_{B}$ the topology on $\SC$ generated by the family of seminorms
\begin{equation}\label{eq:seminormTB}
p_{I,j}(f)=\sup_{v\in C(I,\B_j)}\int_I\big\|\tilde f(s,v(s))\big\|\,ds\,,\quad f\in\SC\, ,
\end{equation}
where $I=[q_1,q_2]$, $q_1,q_2\in\Q$  and $j\in\N$. One has that
$\left(\SC,\T_{B}\right)$ is a locally convex metric space.
\end{defn}

\begin{defn}[Topologies $\T_{D}$ and $\sigma_{D}$]\label{def:TUD}
Let $D$ be a countable and dense subset of $\R^N$. We call $\T_{D}$ (resp.~$\sigma_{D}$) the topology on $\SC$  generated by the family of seminorms
\begin{equation*}
p_{I, y}(f)=\int_I \big\|\tilde f(s,\tilde y)\big\|\,ds\;\;\left(\text{resp.~}\tilde p_{I, y}(f)=\left\|\,\int_I \tilde f(s,\tilde y)\, ds\,\right\|=\sup_{x\in\oU} \left|\,\int_I f(s,x,y)\, ds\,\right| \right)
\end{equation*}
for each $f\in\SC$, $ y\in D$, and $I=[q_1,q_2],\,  q_1,q_2\in\Q $, where $\tilde y$ stands for the identically equal to $y$ map defined on $\oU$.  One has that
$\left(\SC,\T_{D}\right)$ and $\left(\SC,\sigma_{D}\right)$ are locally convex metric spaces.
\end{defn}
Furthermore, in the present PDEs context it makes sense to consider the following so-called {\em pointwise topologies\/}, specially useful in applications.
\begin{defn}[Topologies $\T_{\widetilde DD}$ and $\sigma_{\widetilde DD}$]\label{def:TtildeDD}
Let $\widetilde D$ be a countable and dense subset of $U\subset \R^M$ and let $D$ be a countable and dense subset of $\R^N$. We call $\T_{\widetilde DD}$ (resp.~$\sigma_{\widetilde DD}$) the topology on $\SC$  generated by the family of seminorms
\begin{equation*}
p_{I, x,y}(f)= \int_I |f(s,x,y)|\, ds\, \;\;\left(\text{resp.~}\tilde p_{I, x,y}(f)= \left|\,\int_I f(s,x,y)\, ds\,\right|\right)
\end{equation*}
for each $f\in\SC$, $x\in\widetilde D$,  $ y\in D$, and $I=[q_1,q_2],\,  q_1,q_2\in\Q $.  One has that $\left(\SC,\T_{\widetilde DD}\right)$ and
$\left(\SC,\sigma_{\widetilde DD}\right)$ are locally convex metric spaces.
\end{defn}

Finally, we consider two more metric topologies based on suitable sets of moduli of continuity, trying to adapt the analogous classes proposed and studied  in~\cite{paper:LNO1} and~\cite{paper:LNO2} for ODEs. However, things are more complicated in this PDEs setting and we need to combine two families of moduli of continuity, as in the next definition. The reasons for this precise construction will become clear after the proof of Proposition~\ref{prop:mild solution in K}.

\begin{defn}[Suitable set of moduli of continuity]\label{def:ssmc}
We call  \emph{a suitable set of moduli of continuity}  every $\Theta = (\Theta_1,\Theta_2)$ formed by two countable  sets of non-decreasing continuous functions,
\begin{equation*}
\Theta_1=\left\{\theta^{I}_j \in C(\R^+\!, \R^+)\mid j\in\N, \ I=[q_1,q_2], \ q_1,q_2\in\Q  \right\}
\end{equation*}
such that $\theta^{I}_j(0)=0$ for every $\theta^{I}_j\in\Theta_1$, and  the relation of partial order given~by
\begin{equation*}
\theta^{I_1}_{j_1}\le\theta^{I_2}_{j_2}\quad \text{whenever } I_1\subseteq I_2\,  \text{ and } j_1\le j_2 \,
\end{equation*}
holds, and
\begin{equation*}
\Theta_2=\left\{\theta^{I,\nu}_j \in C(\R^+\!, \R^+) \,\Big|\; j,\nu\in\N, \ I=[q_1,q_2], \ q_1,q_2\in\Q, \ \nu\geq 1 \;\text{s.t.}\; q_1+\frac{1}{\nu}<q_2 \right\}
\end{equation*}
such that $\theta^{I,\nu}_j(0)=0$ for every $\theta^{I,\nu}_j\in \Theta_2$ and
\begin{equation}\label{def:modCont}
\theta^{I_1,\nu_1}_{j_1}\le\theta^{I_2,\nu_2}_{j_2}\quad \text{whenever } I_1\subseteq I_2\,,\; \nu_1\le \nu_2\,, \text{ and } j_1\le j_2 \,,
\end{equation}
with the restriction that $1/\nu_1$ is smaller than the length of the interval $I_1$.
\end{defn}
\begin{defn}[Topologies  $\T_{\Theta}$ and $\sigma_{\Theta}$]\label{def:TTheta}
Let $\Theta= (\Theta_1,\Theta_2)$ be a suitable set of moduli of continuity. For each $j\in\N$ and  $I=[q_1,q_2]$ ($ q_1,q_2\in\Q$), let  $\mathcal{H}_j^{I}$ be the set of functions  $v\in C(I,\B_j)$ which satisfy that
for every integer $\nu\geq 1$ such that $q_1+1/\nu<q_2$,
\begin{itemize}
  \item[] $\|v(t_1)-v(t_2)\|\le \theta^{I,\nu}_j(|t_1-t_2|)+\theta_j^I\Big(\displaystyle\frac{1}{\nu}\Big)$ for all $ t_1, t_2\in \Big[q_1+\displaystyle\frac{1}{\nu},q_2\Big]$.
 \end{itemize}
We call $\T_{\Theta}$ (resp.~$\sigma_{\Theta}$)  the topology on $\SC$ generated by the family of seminorms
\begin{equation*}
\begin{split}
p_{I,j}(f)&=\sup_{v\in\mathcal{H}_j^{I}}\int_I\big\|\tilde f(s,v(s))\big\|\,ds\\ \bigg(\text{resp.}\;\, \tilde p_{I,j}(f)&= \sup_{v\in\mathcal{H}_j^{I}}\left\|\,\int_I \tilde f(s,v(s))\,ds\,\right\| =  \sup_{x\in \oU,\,v\in\mathcal{H}_j^{I} }\bigg|\,\int_I  f(s,x,v(s)(x))\,ds\,\bigg|\bigg)
\end{split}
\end{equation*}
for each $ f\in\SC$, $I=[q_1,q_2]$ ($q_1,q_2\in\Q$), and $j\in\N$. One has that $\left(\SC,\T_{\Theta}\right)$ and $\left(\SC,\sigma_{\Theta}\right)$  are locally convex metric spaces.
\end{defn}
Notice that, as well as  $\SC$, also $\LC$ can be endowed with all the previous topologies. Besides,  given any dense and countable sets $\widetilde D\subset U$ and $D\subset \R^N$, and given any suitable  set of moduli of continuity $\Theta$, it is immediate that the following chains of order hold true:
\begin{equation}\label{eq:TopIncl}
\T_{\widetilde DD}\le \T_{D}\le \T_{\Theta}\le
\T_{ B}\quad\text{and}\quad \sigma_{\widetilde DD}\leq \sigma_{D}\le\sigma_{\Theta} \leq \T_{\Theta}\,.
\end{equation}

In the remaining part of this section we are going to give some definitions and properties on the behaviour of the $m$-bounds and/or the $l$-bounds associated to a set $E\subset \LC$. We might consider $m$-bounds just for maps in $\C$, but since the aim is to study the dependance of mild solutions $u(\cdot,f,z)$ of the abstract problems~\eqref{eq:acp} with respect to the map $f\in \LC$ and the initial condition $z\in X$, we restrict ourselves to $E\subset \LC$.

\begin{defn}
We say that:
\begin{itemize}[leftmargin=19pt]
\item[(i)] A set $E\subset\LC$ \emph{has $L^1_{loc}$-bounded} $m$-bounds (resp.~$l$-bounds [possibly with respect to a modulus of continuity $\widehat\theta$ in the variable $x$]), if for every $j\in\N$ there exists a set  $S^j\subset L^1_{loc}$ of $m$-bounds (resp.~$l$-bounds) of the functions of $E$ on $B_j\subset\R^N$, such that $S^j$ is $L^1_{loc}$\nbd-bounded, that is, for every $r>0$,
\begin{equation*}
\sup_{m\in S^j}\int_{-r}^r m\,<\infty\,.
\end{equation*}


\item[(ii)] A set $E\subset\LC$ \emph{has $L^1_{loc}$-equicontinuous} $m$-bounds (resp.~$l$-bounds [possibly with respect to a modulus of continuity $\widehat\theta$ in the variable $x$]), if for every $j\in\N$ there exists a set  $S^j\subset L^1_{loc}$ of $m$-bounds (resp.~$l$-bounds) of the functions of $E$  on $B_j\subset\R^N$, such that $S^j$ is $L^1_{loc}$-equicontinuous, that is, for every $r>0$ and every $\ep>0$, there exists a $\delta=\delta(j,r,\ep)>0$ such that, for all $-r\le t_1\le t_2\le r$ with $t_2-t_1<\delta$,
\begin{equation*}
\sup_{m\in S^j}\int_{t_1}^{t_2} m\,<\ep.
\end{equation*}
\end{itemize}
\label{def:05.07-13:05}
\end{defn}
\begin{rmk}\label{rmk:equicntThenBounded}
According to the previous definitions,  the $L^1_{loc}$-equicontinuity implies the $L^1_{loc}$-boundedness.
\end{rmk}
The next result asserts that the existence of $L^1_{loc}$-bounded or $L^1_{loc}$-equicontinuous $m$-bounds and/or $l$-bounds for a set $E\subset \LC$ is inherited by all the elements in the closure of $E$ with respect to any of the previously introduced topologies. The proof is omitted as it presents only minor differences with respect to the one of Proposition~2.26 in~\cite{Longothesis} (which in turn extends Propositions~4.7 and~4.10 in~\cite{paper:LNO1}).
\begin{prop}\label{prop:propagation-mlbounds}\label{prop:07.07-19:44}
 Let $\T$ be any of the introduced topologies.
\begin{itemize}[leftmargin=18pt]
\item[\rm (i)] If $E\subset \LC$ has $L^1_{loc}$-bounded $m$-bounds  (resp.~$l$-bounds [possibly with respect to a modulus of continuity $\widehat\theta$]) then $\mathrm{cls}_{(\LC,\T)}(E)$ has $L^1_{loc}$-bounded $m$-bounds  (resp.~$l$-bounds).
\item[\rm (ii)]  If $E\subset \LC$  has $L^1_{loc}$-equicontinuous $m$-bounds  (resp.~$l$-bounds [possibly with respect to a modulus of continuity $\widehat\theta$]), then $\mathrm{cls}_{(\LC,\T)}(E)$   has $L^1_{loc}$-equicontinuous $m$-bounds  (resp.~$l$-bounds).
 \end{itemize}
\end{prop}
The fact that a set $E\subset\LC$ has $L^1_{loc}$-bounded $l$-bounds has strong implications on the equivalence of some strong topologies when restricted to $E$. When additional conditions are imposed on $E$, then some weak topologies turn out to be equivalent, too.   We collect the  precise  results in the next theorem.
\begin{thm}\label{thm:equivalencia topologias}
Let $\widetilde D\subset U$ and $D\subset \R^N$ be any  dense and countable subsets, and let $\Theta= (\Theta_1,\Theta_2)$ be any suitable set of moduli of continuity, according to Definition~$\ref{def:ssmc}$. The following assertions hold for a set  $E\subset\LC$:
\begin{itemize}
  \item[\rm{(i)}] If $E$ has $L^1_{loc}$-bounded $l$-bounds, then for any $\T_1,\T_2\in\{\T_{D},\T_{\Theta},\T_{B}\}$,
  \begin{equation*}
(E,\T_1)=(E,\T_2)\qquad\text{and}\qquad\mathrm{cls}_{(\LC,\T_1)}(E)=\mathrm{cls}_{(\LC,\T_2)}(E)\,.
\end{equation*}

  \item[\rm{(ii)}] If $E$ has $L^1_{loc}$-equicontinuous $m$-bounds and $L^1_{loc}$-bounded $l$-bounds, then
\begin{equation*}
(E,\sigma_D)=(E,\sigma_\Theta)\qquad\text{and}\qquad\mathrm{cls}_{(\LC,\sigma_D)}(E)=\mathrm{cls}_{(\LC,\sigma_\Theta)}(E)\,.
\end{equation*}

  \item[\rm{(iii)}] If $E$ has $L^1_{loc}$-bounded $l$-bounds with respect to a modulus of continuity $\widehat \theta$ in the variable $x$, then
\begin{equation*}
(E,\sigma_{\widetilde DD})=(E,\sigma_D)\qquad\text{and}\qquad\mathrm{cls}_{(\LC,\sigma_{\widetilde DD})}(E)=\mathrm{cls}_{(\LC,\sigma_D)}(E)\,,
\end{equation*}
and, for any $\T\in\{\T_{D},\T_{\Theta},\T_{B}\}$,
\begin{equation*}
(E,\T_{\widetilde DD})=(E,\T)\qquad\text{and}\qquad\mathrm{cls}_{(\LC,\T_{\widetilde DD})}(E)=\mathrm{cls}_{(\LC,\T)}(E)\,.
\end{equation*}
\end{itemize}

\end{thm}
\begin{proof}
(i) The main ideas come from the proof of Theorem~4.12 in~\cite{paper:LNO1} in an ODEs context, but we include the proof for the sake of completeness.  Due to~\eqref{eq:TopIncl}, it suffices to prove that if $(f_n)_{n\geq 1}$ is a sequence of elements of $ E$ converging to some $f$ in $(\LC,\T_{D})$, then $(f_n)_{n\geq 1}$ converges to $f$ in $(\LC,\T_{B})$ too.  Fix a compact interval $I=[q_1,q_2]$ with $q_1,q_2\in\Q$ and fix a $j\in\N$. For each ${n\geq 1}$, let $l_n^j\in L^1_{loc}$ be the optimal $l$-bound of $f_n$ on $B_j$, and $l_f^j\in L^1_{loc}$ be the one for $f$,  defined as in~\eqref{eqOptimalMLbound}. By hypothesis and by Proposition~\ref{prop:propagation-mlbounds}, there exists a $\rho>0$ such that
\begin{equation}\label{eq:bounded l bounds rho}
\sup_{{n\geq 1}}\int_I l^j_n(s)\, ds<\rho<\infty \quad\text{and} \quad \int_I l_f^j(s)\,ds<\rho\, .
\end{equation}
Now, fix $\ep>0$ and consider $\delta=\ep/(3\,\rho)$. Since $ B_j \subset \R^N$ is compact, and $D$ is dense in $\R^N$, there exist  $y_1,\dots, y_{i_0}\in D$ such that $ B_j \subset \bigcup_{i=1}^{i_0} \stackrel{\ \circ} B_\delta\!(y_i)$, where $B_\delta(y)$ denotes the closed ball of $\R^{N}$ of radius $\delta$ centered at $y\in\R^{N}$. Subordinate to the former open cover, there is a continuous partition of unity, that is,  there exist continuous functions $\phi_i:\R^N\to[0,1]$ for $i=1,\dots,{i_0}$   such that
\begin{equation*}
\supp(\phi_i)\subset \stackrel{\ \circ} B_\delta\!(y_i)\, \qquad \mathrm{and}\qquad \sum_{i=1}^{i_0} \phi_i(y)=1 \quad \text{for all }\, y\in B_j \,.
\end{equation*}
Let us define the $\SC$ functions
\begin{equation}\label{eq:f estrella}
f^*_n(t,x,y)=\sum_{i=1}^{i_0}\phi_i(y)\, f_n(t,x,y_i)\quad \text{and}\quad
f^*(t,x,y)=\sum_{i=1}^{i_0} \phi_i(y)\,  f(t,x,y_i)\, .
\end{equation}
Then, for each  $v\in C(I,\B_j)$, note that $v(t)(x)\in B_j$ for all $t\in I$ and $x\in\oU$, and write
\begin{equation*}
\begin{split}
\int_I \big\|\tilde f_{n}(t,v(t)) -   \tilde f(t,v(t))\big\|\,dt &= \int_I \sup_{x\in\oU}\left|f_{n}(t,x,v(t)(x)) -   f(t,x,v(t)(x))\right|\,dt \\&\le\int_I \sup_{x\in\oU}\left|f_{n}(t,x,v(t)(x)) -   f_n^*(t,x,v(t)(x))\right|\,dt\\
			&\hspace{0,5cm}+\int_I \sup_{x\in\oU}\left|f^*_{n}(t,x,v(t)(x)) -   f^*(t,x,v(t)(x))\right|\,dt \\
            &\hspace{0,5cm}+\int_I \sup_{x\in\oU} \left|f^*(t,x,v(t)(x)) -   f(t,x,v(t)(x))\right|\,dt .
\end{split}
\end{equation*}

Let us separately analyze each one  of the three terms in the previous sum. As regards the first one, by the construction of the maps $\phi_i$ we have that for all $n\geq 1$,
\begin{equation}\label{eq:27.10-13:59}
\begin{split}
\int_I &\sup_{x\in \oU}\left|f_{n}(t,x,v(t)(x)) -   f_n^*(t,x,v(t)(x))\right|\,dt
\\&=\int_I \sup_{x\in \oU}\Big|\sum_{i=1}^{i_0} \phi_i(v(t)(x))\left( f_{n}(t,x,v(t)(x))- f_{n}(t,x,y_i)\right)\Big|\,dt\\
&\le \int_I \sup_{x\in \oU} \sum_{i=1}^{i_0} \phi_i(v(t)(x))\left| f_{n}(t,x,v(t)(x))- f_{n}(t,x,y_i)\right|\,dt\\
&\le \int_I \sup_{x\in \oU} \sum_{i=1}^{i_0} \phi_i(v(t)(x))\, l^j_{n}(t) \left| v(t)(x)-y_i\right|\,dt\\
&\le \int_I \sup_{x\in \oU}\sum_{i=1}^{i_0} \phi_i(v(t)(x))\, l^j_{n}(t)\ \delta\,dt= \frac{\ep}{3\,\rho}\int_I l^j_{n}(t)\, dt\le \frac{\ep}{3}\,,
\end{split}
\end{equation}
applying \eqref{eq:bounded l bounds rho} to get the last inequality.
As for the third element in the sum, also the optimal $l$-bound $l_f^j\in L^1_{loc}$ of $f$ on $B_j$ satisfies the integral bound in~\eqref{eq:bounded l bounds rho}.
Therefore, reasoning exactly as in~\eqref{eq:27.10-13:59}, we obtain that
\begin{equation}\label{eq:27.10-14:00}
\int_I \sup_{x\in \oU}\left|f^*(t,x,v(t)(x)) -   f(t,x,v(t)(x))\right|\,dt \le \frac{\ep}{3}\,.
\end{equation}

Finally, for the remaining element in the sum, since $(f_{n})_{n\geq 1}$ converges to $f$ in $\left(\LC,\T_{D}\right)$, there exists an $n_0=n_0(\ep,i_0)$ such that for all $n\geq n_0$,
\[
p_{I,y_i}(f_n-f)=\int_I \sup_{x\in\overline U}\left|f_{n}(t,x,y_i) -   f(t,x,y_i)\right|\,dt<\frac{\ep}{3\,{i_0}}\quad\text{for all }i=1,\dots, {i_0}.
\]
Then, from~\eqref{eq:f estrella} and the fact that $\phi_i(y)\in[0,1]$ for each $y\in\R^N$ we deduce that, if $n\geq n_0$,
\begin{equation}\label{eq:27.10-16:29}
\begin{split}
\int_I &\sup_{x\in\overline U}\left|f^*_{n}(t,x,v(t)(x)) -   f^*(t,x,v(t)(x))\right|\,dt
\\ &\le \sum_{i=1}^{i_0}\int_I \sup_{x\in\overline U}\left|f_{n}(t,x,y_i) -   f(t,x,y_i)\right|\,dt\le \frac{\ep}{3}\,.
\end{split}
\end{equation}
Notice that~\eqref{eq:27.10-13:59}, \eqref{eq:27.10-14:00}, and~\eqref{eq:27.10-16:29} are independent of $v\in C(I,\B_j)$. Thus, gathering~\eqref{eq:27.10-13:59},~\eqref{eq:27.10-14:00} and~\eqref{eq:27.10-16:29}  together, we obtain that  $p_{I,j}(f_n-f)\leq \ep$ for all $n\geq n_0$, for the seminorm $p_{I,j}$ defined in~\eqref{eq:seminormTB}. Thus, $(f_n)_{n\geq 1}$ converges to $f$ in $(\LC,\T_{B})$ and the proof is finished.

(ii) We remark that in this PDEs context we are assuming the $L^1_{loc}$-equicontinuity of the $m$-bounds of $E$, which turns out to be a natural assumption when dealing with weak topologies. Although the main ideas come from the proof of~\cite[Theorem~2.20(ii)]{paper:LNO3} in a context of diffential equations with delay, the technical details are rather different, so that we give a complete proof.

Once more, since the relations~\eqref{eq:TopIncl} hold, it suffices to prove that if $(f_n)_{n\geq 1}\subset E$ converges to some $f$ in $(\LC,\sigma_{D})$, then $(f_n)_{n\geq 1}$ converges to $f$ in $(\LC,\sigma_{\Theta})$ too. To see it, fix a compact interval $I=[q_1,q_2]$ with $q_1,q_2\in\Q$ and fix a $j\in\N$. We want to prove that $\tilde p_{I,j}(f_n-f)\to 0$ as $n\to\infty$ for the seminorm $\tilde p_{I,j}$ given in Definition~\ref{def:TTheta}. Namely,
\begin{equation}\label{eq:semi f_n - f}
\tilde p_{I,j}(f_n-f)=\sup_{x\in \oU,\,v\in\mathcal{H}_j^{I} }\bigg|\int_I  \big(f_n(s,x,v(s)(x))-f(s,x,v(s)(x))\big)\,ds\,\bigg|\,.
\end{equation}
 We keep the notation introduced before. Moreover, for each ${n\geq 1}$, let $m_n^j\in L^1_{loc}$ be the optimal $m$-bound of $f_n$ on $B_j$, and let $m_f^j\in L^1_{loc}$ be the one for $f$, defined as in~\eqref{eqOptimalMLbound}. Recall that the $l$-bounds satisfy relation~\eqref{eq:bounded l bounds rho}, and by Remark~\ref{rmk:equicntThenBounded} and Proposition~\ref{prop:propagation-mlbounds}, there exists a $\rho_0>0$ such that
\begin{equation}\label{eq:bounded m bounds rho 0}
\sup_{{n\geq 1}}\int_I m^j_n(s)\, ds<\rho_0<\infty \quad\text{and} \quad \int_I m_f^j(s)\,ds<\rho_0\, .
\end{equation}

The initial arguments are quite similar to the ones in (i). Fixed $\ep>0$, we consider $\delta=\ep/(3\,\rho)$ for $\rho$ in~\eqref{eq:bounded l bounds rho}, we take a partition of unity just as in (i), and consider the maps in~\eqref{eq:f estrella}. Now, for each  $v\in \mathcal{H}^I_j$ and each $x\in\oU$ we write
\begin{equation}\label{eq:three terms}
\begin{split}
\bigg| &\int_I  \big(f_{n}(t,x,v(t)(x)) -  f(t,x,v(t)(x))\big)\,dt\, \bigg|   \\
&\leq \left| \,\int_I \big(f_{n}(t,x,v(t)(x)) -  f_{n}^*(t,x,v(t)(x))\big)\,dt \,\right|  \\
&	\hspace{0,5cm}+\left| \,\int_I \big(f_{n}^*(t,x,v(t)(x)) -  f^*(t,x,v(t)(x))\big)\,dt\, \right| \\
    &\hspace{0,5cm}         +\left| \,\int_I \big(f^*(t,x,v(t)(x)) -  f(t,x,v(t)(x))\big)\,dt \,\right|\, .
\end{split}
\end{equation}
The first and third terms in the sum are treated similarly as the corresponding ones in (i), using the properties of the partition of unity and bringing the $l$-bounds into play. Both terms turn out to be less than $\ep/3$ for all $n\geq 1$, uniformly for $v\in \mathcal{H}^I_j$ and $x\in\oU$.

The treatment of the second term is more delicate and technical. By~\eqref{eq:f estrella},
\begin{equation}\label{eq:second term}
\begin{split}
\bigg| &\int_I \big(f_{n}^*(t,x,v(t)(x)) -  f^*(t,x,v(t)(x))\big)\,dt\, \bigg|\\
&= \bigg| \int_I \sum_{i=1}^{i_0} \phi_i(v(t)(x))\big( f_{n}(t,x,y_i)- f(t,x,y_i)\big)\,dt\, \bigg|\,.
\end{split}
\end{equation}
The uniform continuity of the functions $\phi_i$ on the compact set $B_j$ helps to deal with this term. Given $\ep^*=\ep/(18\,i_0\,\rho_0)$, for $\rho_0$ in~\eqref{eq:bounded m bounds rho 0}, there exists a $\delta_0>0$ such that, if $y,\tilde y\in B_j$ satisfy $|y-\tilde y|<\delta_0$, then $|\phi_i(y)-\phi_i(\tilde y)|<\ep^*$ for all $i=1,\ldots,i_0$. Then,  we can determine an integer $\nu\geq 1$ large enough so that $q_1+1/\nu<q_2$ and it satisfies the next two conditions:
\begin{itemize}
  \item[($c_1$)] $\theta^I_j\Big(\displaystyle\frac{1}{\nu}\Big)<\displaystyle\frac{\delta_0}{2}$; \vspace{0,2cm}
  \item[($c_2$)] $\displaystyle\int_{q_1}^{q_1+1/\nu} m^j_n(s)\, ds< \displaystyle\frac{\ep}{18\,i_0}\,$ for all $n\geq 1$, and $\,\displaystyle\int_{q_1}^{q_1+1/\nu} m^j_f(s)\, ds< \displaystyle\frac{\ep}{18\,i_0}\,$.
\end{itemize}
Just recall that  $\theta^I_j\in \Theta_1$ (see Definition~\ref{def:ssmc}) is a  modulus of continuity, and we are assuming that $E$ has $L^1_{loc}$-equicontinuous $m$-bounds and Proposition~\ref{prop:propagation-mlbounds} holds. This time considering the modulus of continuity $\theta^{I,\nu}_j\in \Theta_2$, we can find a small enough  $h>0$, with $h\in\Q$ for convenience,  so that
\begin{itemize}
  \item[($c_3$)] $\theta^{I,\nu}_j(h)<\displaystyle\frac{\delta_0}{2}$.
  \end{itemize}
Then, whenever $v\in \mathcal{H}^I_j$ and $t_1,t_2\in [q_1+1/\nu,q_2]$ with $|t_1-t_2|\leq h$, by Definition~\ref{def:TTheta} and conditions $(c_1)$ and $(c_3)$ we have that
  $\|v(t_1)-v(t_2)\|\le \theta^{I,\nu}_j(|t_1-t_2|)+\theta_j^I(1/\nu)\leq \theta^{I,\nu}_j(h)+\theta_j^I(1/\nu)<\delta_0$. By the uniform continuity of the maps $\phi_i$ on $B_j$, for every $x\in \oU$ we have that $|\phi_i(v(t_1)(x))-\phi_i(v(t_2)(x))|<\ep^*$, for all $i=1,\ldots,i_0$.

So, the idea is to split the integral in~\eqref{eq:second term} over the whole interval $I=[q_1,q_2]$ into a collection of integrals over intervals of length $h$ where we can use the previous bound. Note that in any case we have to separate a first small interval, since the previous bounds only work for  $t_1,t_2\in [q_1+1/\nu,q_2]$. This is not a problem, thanks to condition $(c_2)$.

More precisely, to simplify the writing, for the integer $k_0\geq 0$ such that $q_2=q_1+1/\nu+k_0\,h+h_1$, with $0\leq h_1<h$, let $r_k=q_1+1/\nu+k\,h$ for $k=0,\ldots,k_0$, and let $r_{k_0+1}=q_2$. First of all, using that $0\leq \phi_i\leq 1$ and condition $(c_2)$, we have that for all $v\in \mathcal{H}^I_j$, $x\in\oU$, and $n\geq 1$,
\begin{equation*}
\begin{split}
  \bigg| &\int_{q_1}^{q_1+\frac{1}{\nu}} \sum_{i=1}^{i_0} \phi_i(v(t)(x))\big( f_{n}(t,x,y_i)- f(t,x,y_i)\big)\,dt\, \bigg|\\
  &\leq \int_{q_1}^{q_1+\frac{1}{\nu}} \sum_{i=1}^{i_0} \big( |f_{n}(t,x,y_i)|+ |f(t,x,y_i)|\big)\,dt \\
  &\leq   i_0\, \int_{q_1}^{q_1+\frac{1}{\nu}}  \big( m^j_n(t)+ m^j_f(t)\big)\, dt \leq i_0\, \left( \frac{\ep}{18\,i_0}+\frac{\ep}{18\,i_0}\right)=\frac{\ep}{9}\,.
  \end{split}
\end{equation*}
Secondly, we write
\begin{equation*}
\begin{split}
  \bigg| &\int_{q_1+\frac{1}{\nu}}^{q_2} \sum_{i=1}^{i_0} \phi_i(v(t)(x))\big( f_{n}(t,x,y_i)- f(t,x,y_i)\big)\,dt\, \bigg|\\
  &=  \left|\, \sum_{k=0}^{k_0}\int_{r_k}^{r_{k+1}} \sum_{i=1}^{i_0} \big( \phi_i(v(t)(x))-\phi_i(v(r_k)(x)) \big)\,\big( f_{n}(t,x,y_i)- f(t,x,y_i)\big)\,dt \right.\\
  &\hspace{0,5cm}+ \left. \sum_{k=0}^{k_0} \sum_{i=1}^{i_0} \phi_i(v(r_k)(x)) \int_{r_k}^{r_{k+1}}\big( f_{n}(t,x,y_i)- f(t,x,y_i)\big)\,dt \, \right|\\
  &\leq \sum_{k=0}^{k_0}\int_{r_k}^{r_{k+1}} \sum_{i=1}^{i_0} \big| \phi_i(v(t)(x))-\phi_i(v(r_k)(x)) \big|\,\big( m^j_n(t)+m^j_f(t)\big)\,dt
  \\
  &\hspace{0,5cm}+ \sum_{k=0}^{k_0} \sum_{i=1}^{i_0} \phi_i(v(r_k)(x))\, \bigg|\int_{r_k}^{r_{k+1}}\big( f_{n}(t,x,y_i)- f(t,x,y_i)\big)\,dt \,\bigg|
   \\
  &\leq \int_{q_1+\frac{1}{\nu}}^{q_2} i_0\,\ep^*\,\big( m^j_n(t)+m^j_f(t)\big)\,dt+ \sum_{k=0}^{k_0} \sum_{i=1}^{i_0} \left\| \, \int_{r_k}^{r_{k+1}}\big( \tilde f_{n}(t,y_i)- \tilde f(t,y_i)\big)\,dt \,\right\|
  \\
  &\leq i_0\,\ep^*\,2\,\rho_0 + \sum_{k=0}^{k_0} \sum_{i=1}^{i_0} \tilde p_{I_k,y_i}(f_n-f) = \frac{\ep}{9} + \sum_{k=0}^{k_0} \sum_{i=1}^{i_0} \tilde p_{I_k,y_i}(f_n-f)
  \end{split}
\end{equation*}
for the intervals with rational endpoints $I_k=[r_k,r_{k+1}]$ for $k=0,\ldots,k_0$ and the seminorms $\tilde p_{I_k,y_i}$ defining the weak topology $\sigma_D$ (see Definition~\ref{def:TUD}). Since by hypothesis, $f_n\to f$ in this topology, we can find an integer $n_0$ such that for all $n\geq n_0$,  $\tilde p_{I_k,y_i}(f_n-f)\leq \ep/(9\,(k_0+1)\,i_0)$ for each $i=1,\ldots i_0$ and $k=0,\ldots,k_0$.

Therefore, bringing everything together, we obtain that for all $n\geq n_0$, the so-called second term in the sum~\eqref{eq:three terms} is less than or equal to $\ep/3$, and this happens  uniformly for $v\in \mathcal{H}^I_j$ and $x\in\oU$. In all, we have dealt with the three terms in~\eqref{eq:three terms}. Finally, looking at~\eqref{eq:semi f_n - f} we can affirm that $\tilde p_{I,j}(f_n-f)\to 0$ as $n\to\infty$, as we wanted to prove.

(iii)  We give all the details for the equivalence of the weak topologies. By~\eqref{eq:TopIncl}, it suffices to see that $\sigma_D\leq \sigma_{\widetilde DD}$. So, let us assume that $(f_n)_{n\geq 1}\subset E$ satisfies that $f_n\to f\in\LC$ in the topology $\sigma_{\widetilde DD}$ and let us prove that $f_n\to f$ in the topology $\sigma_{D}$ too. Let us fix an interval $I=[q_1,q_2]$ ($q_1,q_2\in \Q$) and a $y_0\in D$. By definition,
\[
\tilde p_{I,y_0}(f_n-f)=\sup_{x\in \oU}\bigg|\int_I \big(f_n(t,x,y_0)-f(t,x,y_0) \big)\,dt\, \bigg|\,.
\]
This time we need to find a partition of unity subordinate to an adequate open cover of $\oU$, so as to be able to work with only a finite collection of points in $\widetilde D$.

Let us choose a $j\geq 1$ so that $y_0\in B_j$. Once again, we keep the notation introduced before for the optimal $l$-bounds, $l_n^j$ for $n\geq 1$ and $l_f^j$, but note that this time we are assuming that the maps in $E$ are Lipschitz Carath\'{e}odory with respect to a modulus of continuity $\widehat \theta$ in the variable $x$ (see~\ref{LT}) and the optimal $l$-bounds are defined as in~\eqref{eq:Lbound-theta}. We can assume that~\eqref{eq:bounded l bounds rho} holds.

Now, given $\ep>0$, we can find a $\delta>0$ so that $\widehat \theta(\delta)<\ep / (3\,\rho)$, for $\rho$ the one in~\eqref{eq:bounded l bounds rho}. For this $\delta$, since $\widetilde D$ is dense in $U$, and $\oU$ is compact, there exist $x_1,\ldots,x_{k}\in \widetilde D$ such that $\oU\subset\bigcup_{i=1}^{k} \stackrel{\ \circ} B_\delta\!\!(x_i)$, where $B_\delta(x)$ denotes the closed ball of $\R^{M}$ of radius $\delta$ centered at $x\in\R^{M}$. Subordinate to this open cover, there exists a continuous partition of unity, that is,  continuous functions $\varphi_i:\R^M\to[0,1]$ for $i=1,\dots,{k}$  such that
\begin{equation*}
\supp(\varphi_i)\subset \stackrel{\ \circ} B_\delta\!(x_i)\, \qquad \mathrm{and}\qquad \sum_{i=1}^{k} \varphi_i(x)=1 \quad \text{for all }\, x\in \oU .
\end{equation*}
Here we define the $\LC$ functions
\begin{equation*}
f^*_n(t,x,y)=\sum_{i=1}^{k}\varphi_i(x)\, f_n(t,x_i,y)\quad \text{and}\quad
f^*(t,x,y)=\sum_{i=1}^{k} \varphi_i(x)\,  f(t,x_i,y)\, .
\end{equation*}
Then, for each $x\in\oU$ we write
\begin{equation*}
\begin{split}
\bigg| &\int_I  \big(f_{n}(t,x,y_0) -  f(t,x,y_0)\big)\,dt\, \bigg|   \leq \left| \,\int_I \big(f_{n}(t,x,y_0) -  f_{n}^*(t,x,y_0)\big)\,dt \,\right|  \\
&			+\left| \,\int_I \big(f_{n}^*(t,x,y_0) -  f^*(t,x,y_0)\big)\,dt\, \right| +\left| \,\int_I \big(f^*(t,x,y_0) -  f(t,x,y_0)\big)\,dt \,\right|\, .
\end{split}
\end{equation*}
We have a sum of three terms, and the first and the third ones are treated in a similar fashion. We just write down the details for the first term. For each $n\geq 1$, by~\ref{LT}, \eqref{eq:bounded l bounds rho}, the properties of the partition of unity, and the choice of $\delta$, we have that for every $x\in \oU$,
\begin{equation*}
\begin{split}
\bigg| \int_I &\big(f_{n}(t,x,y_0) -  f_{n}^*(t,x,y_0)\big)\,dt \,\bigg| =
\bigg| \int_I  \sum_{i=1}^{k}\varphi_i(x)\,\big( f_n(t,x,y_0)- f_n(t,x_i,y_0) \big)\,dt \,\bigg| \\
&	\leq \int_I  \sum_{i=1}^{k}\varphi_i(x)\, l_n^j(t)\,\widehat \theta(|x-x_i|)\,dt \leq \, \widehat \theta(\delta)\int_I  l_n^j(t)\,dt\leq \frac{\ep}{3}	\, .
\end{split}
\end{equation*}
As for the second term, for every $x\in \oU$,
\begin{equation*}
\begin{split}
\bigg| \int_I &\big(f_{n}^*(t,x,y_0) -  f^*(t,x,y_0)\big)\,dt\, \bigg| =
\bigg| \int_I  \sum_{i=1}^{k}\varphi_i(x)\,\big( f_n(t,x_i,y_0)- f(t,x_i,y_0) \big)\,dt \,\bigg| \\
&	=  \bigg| \sum_{i=1}^{k}\varphi_i(x) \int_I \big( f_n(t,x_i,y_0)- f(t,x_i,y_0) \big) \,dt\, \bigg|
\leq \sum_{i=1}^{k} \tilde p_{I,x_i,y_0}(f_n-f)	
\end{split}
\end{equation*}
for the seminorms $\tilde p_{I,x_i,y_0}$ generating the weak topology $\sigma_{\widetilde DD}$ (see Definition~\ref{def:TtildeDD}).
Since $f_n\to f$ in this topology, there exists an $n_0$ such that for all $n\geq n_0$, $\tilde p_{I,x_i,y_0}(f_n-f)\leq \ep/(3\,k)$ for all $i=1,\ldots,k$. As a consequence, we have that $\tilde p_{I,y_0}(f_n-f)\leq \ep$ for every $n\geq n_0$, and we are done. 	

As for the strong topologies, note that by (i), it suffices to prove the equivalence of $\T_{\widetilde DD}$ and $\T_{D}$. We skip the details since only minor modifications are needed over the proof just done for the weak topologies. The proof  is finished.
\end{proof}

\section{Topologies for the continuous variation of mild solutions}\label{secContMildSolutions}
In this section we prove the continuity of the mild solutions $u(\cdot,f,z)$ of the abstract problems~\eqref{eq:acp} associated to the parabolic problems~\eqref{eq:pde}, with respect to the map $f\in \LC$, under certain topologies, and the initial condition $z\in X$, for $X=C(\oU,\R^N)$ or $X=C_0(\oU,\R^N)$ depending on the boundary conditions (see Remark~\ref{rmk:notation X}).

We will carry out the study, first, for a class $E\subset \LC$ with $L^1_{loc}$-equicontinuous $m$-bounds, introducing a new  topology $\sigma_{\Theta\mathcal{R}}$  in Definition~\ref{def:TThetaR} for an appropriate choice of $\Theta$ and $\mathcal{R}$ related to the set $E$, in an attempt to have a topology of continuity as coarse as possible, as well as to recover, in a precise sense explained later, some properties of compactness which held in the ODEs case.  And, second, for a class $E\subset \LC$ with $L^1_{loc}$-bounded $l$-bounds, with respect to any strong topology $\T_D$, $\T_\Theta$ or  $\T_B$, all of them equivalent in this case according to Theorem~\ref{thm:equivalencia topologias}(i).

In any of the two cases, we will denote   by $m_f^j\in L^1_{loc}$ (resp.\ $l_f^j\in L^1_{loc}$) the optimal $m$-bound (resp.\ $l$-bound) of $f$ on  the compact ball $B_j\subset \R^N$. Generally, we will just write $m_n^j$ or $l_n^j$ for the optimal $m$-bound or $l$-bound of a map $f_n$.
\par\medskip\noindent
\subsection{Topologies of continuity when $E$ has $L^1_{loc}$-equicontinuous $m$-bounds}

Let us start by considering $E \subset \LC$ with $L^1_{loc}$-equicontinuous $m$-bounds. Associated to the family $E$, in Proposition~\ref{prop:mild solution in K} we determine a suitable set of moduli of continuity  $\Theta$ as in  Definition~\ref{def:ssmc} and a suitable set of radii (see Definition~\ref{def:ssor}) so that a new topology $\T_{\Theta\mathcal{R}}$, as well as its weak version $\sigma_{\Theta\mathcal{R}}$, are built, and they are both topologies of continuity.
\begin{defn}[Suitable set of radii]\label{def:ssor}
We call  \emph{a suitable set of radii}  every countable set of positive constants $\mathcal{R}=\{R_j^I \mid j\in\N, \ I=[q_1,q_2], \ q_1,q_2\in\Q \}$ such that $R_{j_1}^{I_1}\leq R_{j_2}^{I_2}$ whenever $I_1\subseteq I_2$ and $ j_1\le j_2$.
\end{defn}

Let us now introduce a new family of sets  $\K_j^{I}$ which satisfy $\K_j^{I}\subset \mathcal{H}_j^{I}$ for the sets $\mathcal{H}_j^{I}$ given in Definition~\ref{def:TTheta}. As in the previous section, we denote $\B_j=C(\oU,B_j)$ and we will write $\B_j^X$ to mean $\B_j\cap X$. Also, ${\rm d}$ denotes the usual distance between a map $z$  and  a set of maps $F$ in $C(\oU,\R^N)$, ${\rm d}(z,F):=\inf_{\tilde z\in F} \|z-\tilde z\|$.

\begin{defn}\label{def:sets K}
Let $\Theta= (\Theta_1,\Theta_2)$ be a suitable set of moduli of continuity and let $\mathcal{R}$ be a suitable set of radii. For each $j\in\N$ and  $I=[q_1,q_2]$ ($ q_1,q_2\in\Q$), we define $\K_j^{I}$ as the set of functions  $v\in C(I,\B_j)$ which satisfy that, for every interval $J=[r_1,r_2]\subseteq I$ ($r_1,r_2\in\Q$) and for every integer $\nu\geq 1$ such that $r_1+1/\nu<r_2$,
\begin{itemize}
  \item[\rm{(i)}] $\|v(t_1)-v(t_2)\|\le \theta^{J,\nu}_j(|t_1-t_2|)+\theta_j^J\Big(\displaystyle\frac{1}{\nu}\Big)$ for all $ t_1, t_2\in \Big[r_1+\displaystyle\frac{1}{\nu},r_2\Big]$;
  \item[\rm{(ii)}] ${\rm d}\big(v(t), e^{1/\nu A}\,\B_{R_j^J}^X\big)\leq \theta_j^J\Big(\displaystyle\frac{1}{\nu}\Big)$ for all $t\in \Big[r_1+\displaystyle\frac{1}{\nu},r_2\Big]$.
\end{itemize}
\end{defn}
\begin{rmks}\label{rmk:intervalo menor}
1. With the previous definition it is easy to check that,  if $v\in \K_j^{I}$ for some $j\in\N$ and  $I=[q_1,q_2]$ ($q_1,q_2\in\Q$), and $J=[r_1,r_2]\subset I$ ($r_1,r_2\in\Q$), then the restriction of $v$ to $J$ satisfies $v|_J\in \K_j^{J}$.

2. One may wonder whether the sets $\K_j^{I}$ are nonempty, since the construction is rather restrictive. For the specific choice of $\Theta$ and $\mathcal{R}$ that will be made in the next result, these sets are nonempty since mild solutions lie therein. To some extent, this is a generalization of what happened in the ODEs context, where the moduli of continuity associated to a set $E$ with $L^1_{loc}$-equicontinuous $m$-bounds (see~\cite[Definition~5.3]{paper:LNO1}) were the appropriate ones for the solutions. As we are going to see, in the present PDEs context some properties of the semigroup of operators $(e^{tA})_{t\geq 0}$ come into play in a natural way, and this is the main reason for the combination of a pair of sets of moduli of continuity. Basically, the pair $\Theta=(\Theta_1,\Theta_2)$ depends upon the $L^1_{loc}$-equicontinuity of the $m$-bounds, but also upon the boundedness of $(e^{tA})_{t\geq 0}$  and  the uniform continuity of $e^{tA}$ for $t$ in compact intervals away from $0$. This justifies the presence of the additional parameter $\nu$ in the set $\Theta_2$.

3.  Whereas the sets $\K_j^{I}$ in the ODEs case were compact, here we can only guarantee that given $I=[q_1,q_2]$ and $j\geq 1$, the restriction of the set of maps $\K_j^{I}$ to each interval $[q_1+\delta,q_2]$ is compact for every $\delta >0$, i.e., the set $\{v|_{[q_1+\delta,q_2]}\mid v \in \K_j^{I}\}$ is compact. To see it, note that it is closed and apply Arzel\`{a}-Ascoli's theorem. From (i) we can deduce the uniform equicontinuity of this set of maps on $[q_1+\delta,q_2]$ by choosing $\nu$ large enough,  and from (ii), since the operators $e^{1/\nu A}$ are compact, we can deduce that for each fixed $t\in [q_1+\delta,q_2]$ the image $v(t)$ of each $v \in \K_j^{I}$  lies within an arbitrarily small distance from the  set $e^{1/\nu A}\,\B_{R_j^J}^X$  of $X$,  which is precompact, once more by taking $\nu$ as big as needed.
\end{rmks}

\begin{prop}\label{prop:mild solution in K}
Let $E\subset \LC$ have $L^1_{loc}$-equicontinuous $m$-bounds. Then, there is an associated suitable set  of moduli of continuity $\Theta=(\Theta_1,\Theta_2)$, with $\Theta_1=\big\{\theta^{I}_j \in C(\R^+\!, \R^+)\mid j\in\N, \ I=[q_1,q_2], \ q_1,q_2\in\Q\}$ and $\Theta_2=\big\{\theta^{I,\nu}_j \in C(\R^+\!, \R^+)\mid j,\nu\in\N, \ I=[q_1,q_2], \ q_1,q_2\in\Q, \ \nu\geq 1  \;\text{s.t.}\; q_1+1/\nu<q_2\big\}$, and there is an associated suitable set of radii $\mathcal{R}=\{R_j^I \mid j\in\N, \ I=[q_1,q_2], \ q_1,q_2\in\Q \}$ such that, if  $\K_j^I$ are the sets  given in Definition$\ \ref{def:sets K}$, and  $u(\cdot,f,z)$ is the mild solution of the problem
\begin{equation*}
\left\{\begin{array}{l} u'(t)  =
 A\, u(t)+\tilde f(t,u(t))\,,\quad t>q_1\,,\\
u(q_1)=z\,,
\end{array}\right.
\end{equation*}
on $I=[q_1,q_2]$,  $q_1,q_2\in\Q$, for certain $f\in E$ and $z\in X$, and for some $j\geq 1$, $\|u(t,f,z)\|\leq j$ for all $t\in I$, then $u(\cdot,f,z)\in \K_j^I$.
\end{prop}
\begin{proof}
Let $u(t):=u(t,f,z)$, $t\in I=[q_1,q_2]$. Trying to define a modulus of continuity in $t$, for $t_1, t_2\in I$ with $t_1<t_2$ we write, having in mind~\eqref{eq:variation constants} and Remark~\ref{rmk:puntoPartida}.1,
\begin{equation}\label{eq:equicontinuity}
\begin{split}
u(t_2)-&u(t_1)=e^{(t_2-q_1)A}\,z  +\int_{q_1}^{t_2} e^{(t_2-s)A}\,\tilde f(s,u(s))\,ds - e^{(t_1-q_1)A}\,z \\ & - \int_{q_1}^{t_1} e^{(t_1-s)A}\,\tilde f(s,u(s))\,ds = \big(e^{(t_2-q_1)A}- e^{(t_1-q_1)A}\big)\,z \\
&+ \int_{q_1}^{t_1} \big(e^{(t_2-s)A}-e^{(t_1-s)A}\big)\,\tilde f(s,u(s))\,ds+ \int_{t_1}^{t_2} e^{(t_2-s)A}\,\tilde f(s,u(s))\,ds\,.
\end{split}
\end{equation}
Since the semigroup of operators $(e^{tA})_{t\geq 0}$ is analytic, the map $t\in (0,\infty)\mapsto e^{tA}\in \linear (X)$ is in particular norm continuous. Now, note that, to bound the term $(e^{(t_2-q_1)A}- e^{(t_1-q_1)A})\,z$, we have to keep away from $q_1$. The reason is that, even if the semigroup of operators is strongly continuous, we are considering initial conditions $z$ which vary in a bounded set of $X$, $\|z\|\leq j$, which is not a compact set in $X$, so that we cannot obtain a modulus of continuity at $q_1$ which is valid for all such $z$. What is true is that, writing $\|(e^{(t_2-q_1)A}- e^{(t_1-q_1)A})\,z\|\leq j\,\|e^{(t_2-q_1)A}- e^{(t_1-q_1)A}\|$, given any integer $\nu\geq 1$ such that $q_1+1/\nu<q_2$, we can determine a modulus of continuity for $t_1, t_2\in [q_1+1/\nu,q_2]$ by the uniform continuity of $t\mapsto e^{tA}$ on  compact intervals $[1/\nu,T]$, for each $T>1/\nu$. Namely, the map
\begin{equation}\label{eq:modulo cont exponencial}
\theta^{I,\nu}(s):=\sup_{\substack{t\in [\frac{1}{\nu},q_2-q_1]\\ r\in[0,s]}}  \|e^{tA}- e^{(t+r)A}\|\,,\quad s\geq 0
\end{equation}
is a modulus of continuity and $\|(e^{(t_2-q_1)A}- e^{(t_1-q_1)A})\,z\|\leq j\, \theta^{I,\nu}(t_2-t_1) $ for all $q_1+1/\nu\leq t_1 < t_2\leq q_2$. Besides, it is easy to check that, if $\nu_1\leq \nu_2$ and $I_1\subseteq I_2$ with the restriction that $1/\nu_1$ is smaller than the length of the interval $I_1$, then  $\theta^{I_1,\nu_1}\leq \theta^{I_2,\nu_2}$.
 This is independent of $f\in E$ and $z\in X$ with $ \|z\|\leq j$.

Now we tackle the integral terms. Taking relations~\eqref{eq:bound semigroup} and~\eqref{eq:bound} into account, the last term in~\eqref{eq:equicontinuity} is bounded for each $q_1\leq t_1<t_2\leq q_2$ as follows,
\begin{align*}
\left\|\,\int_{t_1}^{t_2} e^{(t_2-s)A}\,\tilde f(s,u(s))\,ds\,\right\|&\leq \int_{t_1}^{t_2} \big\|e^{(t_2-s)A}\big\|\,\big\|\tilde f(s,u(s))\big\|\,ds\\
&\leq N\int_{t_1}^{t_2} m_f^j(s)\,ds\leq  \theta_j^I (t_2-t_1)\,,
\end{align*}
where $\theta_j^I$ is  a slight modification of the   modulus of continuity associated to $E$ given in Definition~5.3 in~\cite{paper:LNO1} in an ODEs context. Precisely, here we take for convenience
\begin{equation}\label{eq:thetaODEs}
\theta_j^I(s):=2\,N\sup_{t\in I, g\in E} \int_t^{t+s} m_g^j\,, \quad s\geq 0\,.
\end{equation}
The fact that this is a modulus of continuity  is due to the hypothesis that the family $E$ has $L^1_{loc}$-equicontinuous $m$-bounds. This is again valid for all $f\in E$ and $z\in X$ with $ \|z\|\leq j$, as far as the solution $u(t)$ keeps bounded by $j$ on the interval $I$.

Finally, we deal with the remaining integral term in \eqref{eq:equicontinuity}.   We split the integral in two. For each $t_1 < t_2$ in $[q_1+1/\nu,q_2]$ we write
\begin{equation*}
\left\|\,\int_{t_1-\frac{1}{\nu}}^{t_1} \big(e^{(t_2-s)A}-e^{(t_1-s)A}\big)\,\tilde f(s,u(s))\,ds\,\right\|\leq 2\,N \int_{t_1-\frac{1}{\nu}}^{t_1}  m_f^j(s)\,ds\leq \theta_j^I\Big(\frac{1}{\nu}\Big)\,.
\end{equation*}
Finally, whenever $t_1>q_1+1/\nu$, we bound the term
\begin{equation}\label{eq:term}
\begin{split}
&\hspace{-0.5cm}\left\| \, \int_{q_1}^{t_1-\frac{1}{\nu}} \big(e^{(t_2-s)A}-e^{(t_1-s)A}\big)\,\tilde f(s,u(s))\,ds \,\right\|\\
&\leq \sup_{s\in [q_1,t_1-\frac{1}{\nu}]} \big\|e^{(t_2-s)A}-e^{(t_1-s)A}\big\|\,  \sup_{g\in E}\int_I m^j_g \leq  \theta^{I,\nu}(t_2-t_1)\,\sup_{g\in E}\int_I m^j_g  \, ,
\end{split}
\end{equation}
for $\theta^{I,\nu}$ given in \eqref{eq:modulo cont exponencial}. Recall that $E$ has $L^1_{loc}$-bounded $m$-bounds by  Remark~\ref{rmk:equicntThenBounded}. Bringing everything together, we can  conclude that appropriate moduli of continuity $\theta^{I,\nu}_j \in C(\R^+\!, \R^+)$ for each  $j\in\N, \  I=[q_1,q_2], \ q_1,q_2\in\Q$, and $\nu\geq 1$ such that $q_1+1/\nu<q_2$ can be determined, so that the order condition~\eqref{def:modCont} holds and $\|u(t_1)-u(t_2)\|\le \theta^{I,\nu}_j(|t_1-t_2|)+\theta^{I}_j(1/\nu)$ for all $t_1,t_2\in [q_1+1/\nu,q_2]$. More specifically,
\[
\theta^{I,\nu}_j (s):= j\, \theta^{I,\nu}(s) + \theta_j^I(s) + \theta^{I,\nu}(s)\,\sup_{g\in E}\int_I m^j_g\,,\quad s\geq 0\,.
\]
The previous arguments work identically if we consider any other interval $J=[r_1,r_2]\subset I$ ($r_1,r_2\in\Q$) and every integer $\nu\geq 1$ such that $r_1+1/\nu<r_2$, by recalling Remark~\ref{rmk:puntoPartida}.2 and writing
$u(t)=e^{(t-r_1)A}\,u(r_1) +\int_{r_1}^t e^{(t-s)A}\,\tilde f(s,u(s))\,ds$, for $t\in[r_1,r_2]$.
Therefore,  (i) in Definition~\ref{def:sets K} holds for $u(t)$.

As for (ii), we have to determine a suitable set of radii $\mathcal{R}$ according to Definition~\ref{def:ssor}, so that for every interval $J=[r_1,r_2]\subseteq I$ and for every $\nu\geq 1$ such that $r_1+1/\nu<r_2$ it holds that ${\rm d}\big(u(t), e^{1/\nu A}\,\B_{R_j^J}^X\big)\leq \theta_j^J(1/\nu)$ for all $t\in [r_1+1/\nu,r_2]$. As before, it suffices to argue for $I=[q_1,q_2]$ and $\nu\geq 1$ such that $q_1+1/\nu<q_2$.

Since $e^{(t+s)A}=e^{tA}\,e^{sA}$ for all $t,s\geq 0$, and $e^{1/\nu A}\in \mathcal{L}(X)$, for each $t \in [q_1+1/\nu,q_2]$ we can write
\begin{align*}
u(t)&=e^{\frac{1}{\nu} A}\,\left(e^{(t-q_1-\frac{1}{\nu})A}\,z +\int_{q_1}^{t-\frac{1}{\nu}}  e^{(t-s-\frac{1}{\nu})A}\,\tilde f(s,u(s))\,ds\right)
\\ &\hspace{0,4cm}+ \int_{t-\frac{1}{\nu}}^t e^{(t-s)A}\,\tilde f(s,u(s))\,ds\,.
\end{align*}
Now, by \eqref{eq:bound semigroup} and~\eqref{eq:bound} we can bound
\[
\left\| \,e^{(t-q_1-\frac{1}{\nu})A}\,z +\int_{q_1}^{t-\frac{1}{\nu}}  e^{(t-s-\frac{1}{\nu})A}\,\tilde f(s,u(s))\,ds\,\right\|\leq
N\,j +N\,\sup_{g\in E}\int_I m^j_g=:R_j^I\,.
\]
Before we proceed, note that the constants $R_j^I>0$ defined in this way satisfy that $R_{j_1}^{I_1}\leq R_{j_2}^{I_2}$ whenever $I_1\subseteq I_2$ and $ j_1\le j_2$. Thus, the set $\mathcal{R}=\{R_j^I \mid j\in\N, \ I=[q_1,q_2], \ q_1,q_2\in\Q \}$ is a suitable set of radii associated to the family $E$.

Finally, it is immediate  to conclude that
\[
{\rm d}\big(u(t), e^{\frac{1}{\nu} A}\,\B_{R_j^I}^X\big)\leq \left\|\,\int_{t-\frac{1}{\nu}}^t e^{(t-s)A}\,\tilde f(s,u(s))\,ds\,\right\|\leq\theta_j^I\Big(\frac{1}{\nu}\Big),
\]
for the modulus in~\eqref{eq:thetaODEs}, as we wanted. Therefore, $u(\cdot,f,z)\in \K_j^I$ and the proof is finished.
\end{proof}
\begin{rmk}
From now on, we will refer to $\Theta=(\Theta_1,\Theta_2)$ and $\mathcal{R}$ as {\em the associated set of moduli of continuity and set of radii\/}, respectively, for a family $E\subset \LC$ with $L^1_{loc}$-equicontinuous $m$-bounds.
\end{rmk}
We can now define the new topologies $\T_{\Theta\mathcal{R}}$ and $\sigma_{\Theta\mathcal{R}}$.
\begin{defn}[Topologies  $\T_{\Theta\mathcal{R}}$ and $\sigma_{\Theta\mathcal{R}}$]\label{def:TThetaR}
Let $E\subset \LC$ have $L^1_{loc}$-equicontinuous $m$-bounds and let $\Theta=(\Theta_1,\Theta_2)$ and $\mathcal{R}$ be the associated sets of moduli of continuity and radii, respectively.  For each $j\in\N$ and  $I=[q_1,q_2]$ ($ q_1,q_2\in\Q$), let  $\K_j^{I}$ be the set of functions defined in Definition~\ref{def:sets K}.
We call  $\T_{\Theta\mathcal{R}}$ (resp.~$\sigma_{\Theta\mathcal{R}}$)  the topology on $\SC$ generated by the family of seminorms
\begin{equation*}
\begin{split}
p_{I,j}(f)&=\sup_{v\in\K_j^{I}}\int_I\big\|\tilde f(s,v(s))\big\|\,ds\\ \bigg(\text{resp.}\;\, \tilde p_{I,j}(f)&= \sup_{v\in\K_j^{I}}\left\|\,\int_I \tilde f(s,v(s))\,ds\,\right\| =  \sup_{x\in \oU,\,v\in\K_j^{I} }\bigg|\,\int_I  f(s,x,v(s)(x))\,ds\,\bigg|\bigg)
\end{split}
\end{equation*}
for each $ f\in\SC$, $I=[q_1,q_2]$ ($q_1,q_2\in\Q$), and $j\in\N$.
\end{defn}
One can prove that $\left(\SC,\T_{\Theta\mathcal{R}}\right)$ and $\left(\SC,\sigma_{\Theta\mathcal{R}}\right)$  are locally convex metric spaces.

\begin{lem}\label{lem:conv intervalo menor}
Let $E\subset \LC$ have $L^1_{loc}$-equicontinuous $m$-bounds, let $\Theta=(\Theta_1,\Theta_2)$ and $\mathcal{R}$ be the associated sets of moduli of continuity and radii, respectively, and consider the associated topology  $\sigma_{\Theta\mathcal{R}}$. Let $(f_n)_{n\geq 1}\subset \SC$ be such that $f_n\to f$ in $(\SC,\sigma_{\Theta\mathcal{R}})$ as $n\to\infty$, and take $j\geq 1$,  $I=[q_1,q_2]$ with $q_1,q_2\in\Q$. Then, for all $r_1, r_2\in\Q$ with $q_1\leq r_1\leq r_2\leq q_2$,
\[
\lim_{n\to \infty} \sup_{v\in \K_{j}^{I}}\left\|\, \int_{r_1}^{r_2} \big( \tilde f_n(s,v(s)) -\tilde f(s,v(s)) \big)\,ds\,\right\|=0\,.
\]
\end{lem}
\begin{proof}
Let $J= [r_1,r_2]\subset I$. It suffices to note that, by Remark~\ref{rmk:intervalo menor}.1, for each $v\in \K_{j}^{I}$,
\[
\left\|\,\int_{r_1}^{r_2} \big( \tilde f_n(s, v(s)) -\tilde f(s, v(s)) \big)\,ds\,\right\|\leq \tilde p_{J,j}(f_n-f)\to 0 \quad \text{ as }\; n\to\infty\,.
\]
\end{proof}
We are in a position to state the first continuity result with respect to the weak topology $\sigma_{\Theta\mathcal{R}}$. The main arguments in the proof are the same as the ones in Theorem~3.8(i) in~\cite{paper:LNO2}, but technical difficulties arise due to the presence of the semigroup of operators $(e^{tA})_{t\geq 0}$. Recall that $I_{f,z}$ denotes the maximal interval of definition of the mild solution $u(t,f,z)$ of problem~\eqref{eq:acp}. We refer the reader to Section~\ref{secCaratheodoryPDEs} for more details.
\begin{thm}\label{thm:CONTIN equic m bounds}
Let $E\subset \LC$ have  $L^1_{loc}$-equicontinuous $m$-bounds, associated suitable set of moduli of continuity $\Theta$, and associated suitable set of radii $\mathcal{R}$.  Assume that $(f_n)_{n\geq 1}\subset E$ converges to  $f$ in $(\LC,\sigma_{\Theta\mathcal{R}})$ and $(z_n)_{n\geq 1}\subset X$  converges to  $z\in X$, as $n\to\infty$. Then, the mild solutions $u_n(t):=u(t,f_n,z_n)$ converge uniformly as $n\to\infty$ to the mild solution $u(t):=u(t,f,z)$ on every time-interval $[0,T]\subset I_{f,z}$.
\end{thm}
\begin{proof}
For a fixed $T>0$ such that $[0,T]\subset I_{f,z}$, define
\begin{equation}\label{eq:rho}
0<\rho:=1+\max\big\{\,(\|z_n\|)_{n\geq 1},\, \sup\{\|u(t,f,z)\|\mid t\in [0,T]\}\,\big\}\,.
\end{equation}
Now, for each $n\geq 1$ we consider the map $v_n:[0,T]\to X$ defined by
\[
v_n(t):=\left\{\begin{array}{rl} u_n(t)\,, & \text{if} \;\, t\in [0,T_n]\,,\\[0.1cm]
u_n(T_n)\,, & \text{if} \;\,  t\in [T_n,T]\,, \end{array}\right.
\]
where $T_n:=\sup\{t\in [0,T]\mid \|u_n(s)\|\leq \rho \;\text{for all}\; s\in [0,t]\}$. Note that by~\eqref{eq:rho} and the continuity of $u_n$ it  follows that $T_n>0$ for every $n\geq 1$. Fix an integer $j$ such that $j>\rho$ and recall that we denote by $\B_j$ the closed ball of radius $j$ in $C(\oU,\R^N)$.

We affirm that the set $\{v_n(\cdot)\mid n\geq 1 \}$, which is contained in $C([0,T],\B_j)$ by construction, is relatively compact. To see it, by Arzel\`{a}-Ascoli's theorem, we have to check the equicontinuity at any $t\in [0,T]$ as well as the relative compactness  of the set of maps $\{v_n(t)\mid n\geq 1\}\subset C(\oU,\R^N)$ for each fixed $t\in [0,T]$.

In what respects to equicontinuity, for each $n\geq 1$, if $t_1, t_2\in [T_n,T]$, then $v_n(t_2)-v_n(t_1)=0$, so that we only have to worry about $t_1, t_2\in [0,T_n]$. In this case,  for $t_1<t_2$,  similarly to the writing in~\eqref{eq:equicontinuity}, here  we have
\begin{multline}\label{eq:vnt1-vnt2}
v_n(t_2)-v_n(t_1)= \big(e^{t_2 A}- e^{t_1 A}\big)\,z_n \\
+ \int_{0}^{t_1} \big(e^{(t_2-s)A}-e^{(t_1-s)A}\big)\,\tilde f_n(s,v_n(s))\,ds + \int_{t_1}^{t_2} e^{(t_2-s)A}\,\tilde f_n(s,v_n(s))\,ds\,.
\end{multline}
Note that the set $K=\{(z_n)_{n\geq 1},\,z\}$ is a compact set in $X$. With the strong continuity of $(e^{tA})_{t\geq 0}$, we can assert that the map $[0,T]\times K\to X$, $(t,\tilde z)\mapsto e^{tA}\,\tilde z$ is uniformly continuous on $[0,T]\times K$ (see, e.g.,~\cite[Lemma~I.5.2]{book:EN}). This implies that the term $\|\big(e^{t_2 A}- e^{t_1 A}\big)\,z_n\|$ becomes as small as wanted provided that $t_2-t_1$ is small enough, uniformly on $[0,T_n]$ and $n\geq 1$.

About the integral terms in \eqref{eq:vnt1-vnt2}, the last one is treated as the corresponding one in the proof of Proposition~\ref{prop:mild solution in K}, using that $\|e^{tA}\|\leq N$ for all $t\in I=[0,T]$ and that all the maps $f_n\in E$ and $E$ has $L^1_{loc}$-equicontinuous $m$-bounds. That is, fixed an $\ep>0$ we find a  $\delta>0$ such that $\theta^I_j(\delta)<\ep$, for the modulus $\theta^I_j$ in~\eqref{eq:thetaODEs} (which makes sense whether $T\in\Q$ or not) and then,
$\big\|\int_{t_1}^{t_2} e^{(t_2-s)A}\,\tilde f_n(s,v_n(s))\,ds\big\|< \ep$ for all $n\geq 1$, provided that $0<t_2-t_1\leq \delta$. Now, for the first integral term, we take $t_1^* := \max(0,t_1-\delta)$ which equals $0$ if $t_1\leq \delta$, and equals $t_1-\delta$ if $t_1>\delta$, and just as before,  for all $n\geq 1$,
\[
\left\|\,\int_{t_1^*}^{t_1} \big(e^{(t_2-s)A}-e^{(t_1-s)A}\big)\,\tilde f_n(s,v_n(s))\,ds\,\right\|\leq 2\,N \int_{t_1^*}^{t_1} m_n^j(s)  \,ds<\ep\,.
\]
Finally, provided that $t_1>\delta$, we  bound the remaining integral on $[0,t_1-\delta]$ just  as in~\eqref{eq:term}
and we apply the uniform continuity of $e^{tA}$ on $[\delta,T]$ and the fact that $E$ has $L^1_{loc}$-bounded $m$-bounds. In all, we can affirm that we have proved the uniform equicontinuity of the sequence $\{v_n(\cdot)\mid n\geq 1 \}$ on $[0,T]$.


It remains to check the precompactness of the set $\{v_n(t)\mid n\geq 1\}\subset X$ for each fixed $t\in [0,T]$. If $t=0$, then
$\{v_n(0)\mid n\geq 1\}=\{z_n\mid n\geq 1\}$ which by hypothesis is relatively compact. If $t>0$,  it suffices to
 check that given any $\ep>0$ we can determine a precompact set which is within a distance less than $\ep$ from the former set.
But, no matter if $T\in\Q$ or not, by Proposition~\ref{prop:mild solution in K} we can affirm that $(v_n)_{n\geq 1}\subset \K^{[0,T]}_j$, so that in particular the maps $(v_n)_{n\geq 1}$ satisfy (ii) in Definition~\ref{def:sets K} and this is enough for our purposes (see also Remark~\ref{rmk:intervalo menor}.3).

As a consequence, we can assume without loss of generality that $v_n(t)\to v(t)$ uniformly for $t\in [0,T]$ as $n\to\infty$, for a continuous map $v$. Now, define $T_0$ by
\begin{equation}\label{eq:T0}
T_0:=\sup\{t\in [0,T]\mid \|v(s)\|<\rho-1/2  \,\;\text{for all}\,\; s\in[0,t]\}
\end{equation}
and note that $T_0>0$ because of~\eqref{eq:rho}, $\|z_n\|\to \|z\|$ as $n\to\infty$, $v(0)=z$ and $v$ is continuous. Besides, since $v_n$ converges uniformly to $v$ on $[0,T]$, there exists an $n_0\geq 1$ such that $\|v_n(t)\|\leq \rho -1/4$ for all $t\in [0,T_0]$ and $n\geq n_0$. By the definition of the maps $v_n$, this means that $v_n(t)=u_n(t)=u(t,f_n,z_n)$ for all $t\in [0,T_0]$ and  $n\geq n_0$. Therefore, to finish the proof it suffices to check first that $v(t)=u(t,f,z)$ for $t\in [0,T_0]$, and afterwards that $T_0=T$.

To see that $v(t)$ is actually the mild solution for $f$ and $z$ on $[0,T_0]$, write
\begin{equation*}
v_n(t)=e^{tA}\,z_n +\int_0^t e^{(t-s)A}\,\tilde f_n(s,v_n(s))\,ds\,,\quad t\in [0,T_0]\,, \;n\geq n_0\,.
\end{equation*}
We already know that $v_n(t) \to v(t)$ and $e^{tA}\,z_n \to e^{tA}\,z$ as $n\to\infty$, in fact uniformly for $t\in [0,T_0]$. So, if we prove that for each $t\in (0,T_0]$,
\begin{equation}\label{eq:convergence}
\left\|\, \int_0^t e^{(t-s)A}\,\tilde f_n(s,v_n(s))\,ds  - \int_0^t e^{(t-s)A}\,\tilde f(s,v(s))\,ds \, \right\| \to 0 \quad\text{as}\;\, n\to \infty \,,
\end{equation}
we will be done with the first claim. Let us  separately study these two terms:
\[
a_n:=\left\|\, \int_0^t e^{(t-s)A}\,\tilde f_n(s,v_n(s))\,ds  - \int_0^t e^{(t-s)A}\,\tilde f(s,v_n(s))\,ds  \,\right\|
\]
and
\[
b_n:=\left\|\, \int_0^t e^{(t-s)A}\,\tilde f(s,v_n(s))\,ds  - \int_0^t e^{(t-s)A}\,\tilde f(s,v(s))\,ds  \,\right\|\,.
\]
We start with the second one, since it is easier to tackle. Let $l_f^j$ be the optimal $l$-bound of $f$ on  $B_j$ and apply~\eqref{eq:bound semigroup} and Proposition~\ref{prop:f tilde medible}(iv) to write
\begin{equation*}
\begin{split}
b_n &\leq \int_0^t  \big\|e^{(t-s)A}\big\|\,\big\|\tilde f(s,v_n(s))-\tilde f(s,v(s))\big\|  \,ds \\
&\leq  N\int_0^t  l_f^j(s)\,\|v_n(s)  -  v(s)\|\,ds \\
&\leq N \sup_{s\in [0,T_0]} \|v_n(s)  -  v(s)\|\,\int_0^{T_0}\, l_f^j(s)\,ds \to 0 \quad\text{as}\;\, n\to \infty\,,
\end{split}
\end{equation*}
since $v_n$ converges uniformly to $v$ on $[0,T]$.

The first term $a_n$ is more delicate. By Proposition~\ref{prop:propagation-mlbounds}(ii), keeping the agreed notation for the optimal $m$-bounds, we know that fixed an $\ep>0$, there is a $\delta=\delta([0,T_0],j,\ep,N)>0$ such that whenever $t_1, t_2\in [0,T_0]$ satisfy $|t_1-t_2|\leq \delta$,
\[
\int_{t_1}^{t_2} m_n^j\,,\;\, \int_{t_1}^{t_2} m_f^j < \frac{\ep}{2\, N} \quad \text{for all}\;\, n\geq n_0\,,
\]
and besides, we can take a $\rho_0>0$ such that
\[
\int_{0}^{T_0} m_n^j\,,\;\, \int_{0}^{T_0} m_f^j < \rho_0 \quad \text{for all}\;\, n\geq n_0\,.
\]
Now, on the interval $[\delta,T_0]$, $e^{tA}$ is uniformly continuous. Thus, we can take an $h>0$ with the restrictions  $h< \delta$ and $h\in\Q$ for convenience, so that
\[
\big\|e^{t_1 A}-e^{t_2 A}
\big\|< \frac{\ep}{2\, \rho_0} \quad \text{for all}\;\, t_1,\, t_2\in [\delta,T_0] \;\,\text{with}\;\, |t_1-t_2|<h\,.
\]
The first situation we consider is immediate: if $t\in (0,T_0]$ satisfies $t<\delta$, then for all $n\geq n_0$, thanks to the $L^1_{loc}$-equicontinuous $m$-bounds,
\[
a_n\leq \int_0^t
\big\|e^{(t-s)A}
\big\|\,
\big\|\tilde f_n(s,v_n(s))- \tilde f(s,v_n(s))
\big\| \,ds \leq N \left( \int_0^t m_n^j + \int_0^t m_f^j\right) <\ep\,.
\]
If this is not the case, that is, if $t\geq \delta$, then $h$ is the appropriate length of a subinterval where things go well. We write $t-\delta = k_0\, h + h_1$ where $k_0\in\N$ and $0\leq h_1<h$. Then, to get the idea, we split the integral on $[0,t]$ as follows:
\[
\int_0^t \;= \sum_{k=0}^{k_0-1} \int_{kh}^{(k+1)h} \;+ \;\int_{k_0h}^{t-\delta}\;\; +\; \int_{t-\delta}^t
\]
where we understand that the sum is empty if $k_0=0$, so that there are always $k_0$ terms in the sum. Note that, since $t-\delta-k_0\,h = h_1< h<\delta$, the last two integrals can be treated as the one in the first situation just considered ($t<\delta$). For the remaining sum, note that for $s\in [0,k_0\,h]$, $t-s\geq t-k_0\,h=\delta + h_1 \geq \delta$. Then, we do as follows:
\begin{equation*}
\begin{split}
&\left\|\,\sum_{k=0}^{k_0-1} \int_{kh}^{(k+1)h} e^{(t-s)A}\,\big(\tilde f_n(s,v_n(s))-\tilde f(s,v_n(s))\big)\,ds  \,\right\|\\
&= \left\|\,\sum_{k=0}^{k_0-1} \int_{kh}^{(k+1)h} \big(e^{(t-s)A}- e^{(t-kh)A}\big)\,\big(\tilde f_n(s,v_n(s))-\tilde f(s,v_n(s))\big)\,ds \right. \\
&\hspace{0,5cm}+ \left. \sum_{k=0}^{k_0-1}  e^{(t-kh)A} \,\int_{kh}^{(k+1)h} \big(\tilde f_n(s,v_n(s))-\tilde f(s,v_n(s))\big)\,ds  \,\right\|\\
&\leq \sum_{k=0}^{k_0-1}  \int_{kh}^{(k+1)h} \frac{\ep}{2\, \rho_0}\, \big( m_n^j(s) + m_f^j(s)\big)\,ds \\
 &\hspace{0,5cm}+ N \,\sum_{k=0}^{k_0-1}  \left\|\, \int_{kh}^{(k+1)h} \big(\tilde f_n(s,v_n(s))-\tilde f(s,v_n(s))\big)\,ds\,\right\|\\
 &\leq  \frac{\ep}{2 \,\rho_0}\,\int_{0}^{k_0h}  \big( m_n^j(s) + m_f^j(s)\big)\,ds \\
&\hspace{0,5cm} +  N \,\sum_{k=0}^{k_0-1} \sup_{w\in \K_j^I} \left\|\, \int_{kh}^{(k+1)h} \big(\tilde f_n(s,w(s))-\tilde f(s,w(s))\big)\,ds\,\right\|\\
&\leq \ep +  N \,\sum_{k=0}^{k_0-1} \sup_{w\in \K_j^I} \left\|\, \int_{kh}^{(k+1)h} \big(\tilde f_n(s,w(s))-\tilde f(s,w(s))\big)\,ds\,\right\|\,,
\end{split}
\end{equation*}
for the interval $I=[0,k_0\, h]$ with $k_0\,h\in \Q$. Note that once more Proposition~\ref{prop:mild solution in K} justifies that $v_n \in \K_j^I$, for every $n\geq n_0$. At this point we only have to note that for each fixed $t\geq \delta$ we have a finite number of terms in the sum, namely $k_0=k_0(t)$, and each of them tends to $0$ as $n\to\infty$, by  Lemma~\ref{lem:conv intervalo menor}. After this, we can conclude that $v(t)=u(t,f,z)$, the mild solution for $f$ and $z$, on $[0,T_0]$.

To finish the proof, it suffices to check that $T_0=T$. If not, by~\eqref{eq:T0} and the continuity of $v$, we would have that $\|v(T_0)\|=\rho-1/2$, which falls in contradiction with the definition of $\rho$ in~\eqref{eq:rho}. Therefore,  $T_0=T$ and the proof is finished.
\end{proof}
\begin{cor}
Since  $\sigma_{\Theta\mathcal{R}}\le \sigma_{\Theta}$ and $\sigma_{\Theta\mathcal{R}}\le\T_{\Theta\mathcal{R}}\le\T_{\Theta}\le\T_{ B}$, whenever $E\subset \LC$ has $L^1_{loc}$-equicontinuous  $m$-bounds,  as a corollary of the previous theorem we obtain the corresponding continuity result when $(f_n)_{n\geq 1}$ converges to $f$  in any of the topological spaces  $(\LC,\sigma_{\Theta})$, $(\LC,\T_{\Theta\mathcal{R}})$, $(\LC,\T_{\Theta})$ or $(\LC,\T_B)$.
\end{cor}
\par\smallskip\noindent
\subsection{Topologies of continuity when $E$ has $L^1_{loc}$-bounded $l$-bounds}

In this section, we consider  a class $E\subset \LC$ with $L^1_{loc}$-bounded $l$-bounds, and prove the continuity result of mild solutions with respect to $f$ when any strong topology $\T_{D},\T_{\Theta}$ or $\T_{B}$ is considered, all of which are equivalent in this case, as seen in Theorem~\ref{thm:equivalencia topologias}. Before we state the theorem, we include a technical result which will be used in the proof.
\begin{lem}\label{lem:bounded m bounds}
Let $E\subset \LC$ with $L^1_{loc}$-bounded $l$-bounds and let $D$ be any dense and countable subset of $\,\R^N$.  Assume that $(f_n)_{n\geq 1}\subset E$ converges to  $f$ in $(\LC,\T_{D})$ as $n\to\infty$. Then, the family $(f_n)_{n\geq 1}$ has $L^1_{loc}$-bounded $m$-bounds too.
\end{lem}
\begin{proof}
Consider the closed ball $B_j$ in $\R^N$ and take a $y_0\in B_j\cap D$. Then, for each $n\geq 1$, for almost every $t\in \R$, and for all $x\in \oU$ and $y\in B_j$, $|f_n(t,x,y)| \leq |f_n(t,x,y)- f_n(t,x,y_0)| + |f_n(t,x,y_0)| \leq l_n^j(t) \,|y-y_0| + |f_n(t,x,y_0)|\leq 2\,j\,l_n^j(t) +\|\tilde f_n(t,\tilde y_0)\|$. Since the $l$-bounds $(l_n^j)_{n\geq 1}$ are $L^1_{loc}$-bounded, it suffices to check that given any interval $I=[q_1,q_2]$ with $q_1,q_2\in\Q$, $\sup_{n\geq 1}\int_I \|\tilde f_n(t,\tilde y_0)\|\,ds <\infty $. With this purpose, we write $\|\tilde f_n(t,\tilde y_0)\|\leq \|\tilde f_n(t,\tilde y_0)-\tilde f(t,\tilde y_0)\|+\|\tilde f(t,\tilde y_0)\|$. Since $f_n\to f$ in $\T_D$,  $p_{I,y_0}(f_n-f)=\int_I \|\tilde f_n(s,\tilde y_0) - \tilde f(s,\tilde y_0) \|\,ds \to 0$ as $n\to\infty$, and given any $\ep>0$, $p_{I,y_0}(f_n-f)<\ep$ for all $n$ greater than a certain $n_0$. Since a finite number of maps do not matter, the result follows straightaway.
\end{proof}

\begin{thm}\label{thm:CONTIN bounded l bounds}
Let $E\subset \LC$ have $L^1_{loc}$-bounded $l$-bounds and let $D$ be any dense and countable subset of $\,\R^N$.  Assume that $(f_n)_{n\geq 1}\subset E$ converges to $f$ in $(\LC,\T_{D})$ and $(z_n)_{n\geq 1}\subset X$ converges to $z\in X$, as $n\to\infty$. Then, the mild solutions $u_n(t):=u(t,f_n,z_n)$ converge uniformly as $n\to\infty$ to the mild solution $u(t):=u(t,f,z)$ on every time-interval $[0,T]\subset I_{f,z}$.
\end{thm}
\begin{proof}
The main ideas come from the proof of~\cite[Theorem~5.8(i)]{paper:LNO1} in an ODEs context. Nevertheless, there are some technical differences due to the new infinite dimensional scenario. The proof follows the same outline as that of Theorem~\ref{thm:CONTIN equic m bounds}.  Once the maps $v_n:[0,T]\to X$ have been defined on an interval $[0,T]\subset I_{f,z}$, we have to check that they form a relatively compact set in $C([0,T],X)$. For the equicontinuity, whenever $0\leq t_1<t_2\leq T_n$, we have \eqref{eq:vnt1-vnt2} for $v_n(t_2)-v_n(t_1)$.  The term $(e^{t_2A}-e^{t_2A})\,z_n$ therein is treated identically.  Using~\eqref{eq:bound semigroup}, we bound the term
\[
\left\|\,\int_{t_1}^{t_2} e^{(t_2-s)A}\,\tilde f_n(s,v_n(s))\,ds\,\right\|\leq N \int_{t_1}^{t_2} \big\|\tilde f_n(s,v_n(s))\big\|\,ds\,,
\]
 and
\begin{equation}\label{eq:fn}
\int_{t_1}^{t_2} \big\|\tilde f_n(s,v_n(s))\big\|\,ds \leq \int_{t_1}^{t_2} \big\|\tilde f_n(s,v_n(s))-\tilde f(s,v_n(s))\big\|\,ds + \int_{t_1}^{t_2} m_f^j(s)\,ds\,.
\end{equation}
By Theorem~\ref{thm:equivalencia topologias}(i),  $f_n\to f$ as $n\to\infty$ also in the topology $\T_{B}$. Therefore, taking an interval $I=[0,q_2]$ with $T\leq q_2\in\Q$, we have that
\[
p_{I,j}(f_n-f)=\sup_{w\in C(I,\B_j)}\int_I \big\|\tilde f_n(s,w(s))-\tilde f(s,w(s))\big\|\,ds \to 0\quad\text{as}\;\,n\to\infty\,.
\]
Then, fixed any $\varepsilon>0$, and (if necessary) extending  the maps $v_n(t)$ to the interval $[0,q_2]$ just by $v_n(T_n)$ for $t\in [T,q_2]$, it is clear that $v_n\in C(I,\B_j)$ for every $n\geq 1$ and therefore there exists an integer $n_0$ such that for all $n >n_0$, and for all $0\leq t_1<t_2\leq T_n$,
\[
\int_{t_1}^{t_2} \big\|\tilde f_n(s,v_n(s))-\tilde f(s,v_n(s))\big\|\,ds \leq \int_I \big\|\tilde f_n(s,v_n(s))-\tilde f(s,v_n(s))\big\|\,ds <\frac{\ep}{2\,N}\,.
\]
The absolute continuity of the integral on $[0,T]$ for the integrable maps $m_f^j(s)$ and $\|\tilde f_n(s,v_n(s))-\tilde f(s,v_n(s))\|$ for $n=1,\ldots,n_0$  permits us to determine a $\delta>0$ so that the integral of any of these maps over a subinterval in $[0,T]$ of length less than or equal to $\delta$, is less than $\ep/(2\,N)$. In all, we have found a $\delta=\delta(\ep,T)$ such that, if $0<t_2-t_1\leq \delta$, then
$\big\|\int_{t_1}^{t_2} e^{(t_2-s)A}\,\tilde f_n(s,v_n(s))\,ds\big\|< \ep$ for all $n\geq 1$.

Once this $\delta>0$ has been determined, note that $t_1^* := \max(0,t_1-\delta)$ equals $0$ if $t_1\leq \delta$, and equals $t_1-\delta$ if $t_1>\delta$. Then, from the previous argument, we can assert that for all $n\geq 1$,
\[
\left\|\,\int_{t_1^*}^{t_1} \big(e^{(t_2-s)A}-e^{(t_1-s)A}\big)\,\tilde f_n(s,v_n(s))\,ds\,\right\|\leq 2\,N \int_{t_1^*}^{t_1}\big\|\tilde f_n(s,v_n(s))\big\|  \,ds<2\,\ep\,.
\]
Finally, provided that $t_1>\delta$, we still have to treat the  term
\begin{align*}
&\left\|\,\int_{0}^{t_1-\delta} \big(e^{(t_2-s)A}-e^{(t_1-s)A}\big)\,\tilde f_n(s,v_n(s))\,ds\,\right\|
\\ &\hspace{1cm}\leq \sup_{s\in [0,t_1-\delta]}\big\|e^{(t_2-s)A}-e^{(t_1-s)A}\big\| \,\int_{0}^{T} m_n^j (s)\,ds\,.
\end{align*}
Here we apply the uniform continuity of $e^{tA}$ on $[\delta,T]$ and Lemma~\ref{lem:bounded m bounds}. Summing up, we have the uniform equicontinuity on $[0,T]$ for the maps $(v_n)_{n\geq 1}$.

About the relative compactness of the set $\{v_n(t)\mid n\geq 1\}$ for a fixed $t\in [0,T]$, as in the  proof of Theorem~\ref{thm:CONTIN equic m bounds} we note that  for $t=0$ it is known by hypothesis, and if $t>0$ we  can look for a precompact set  within an arbitrarily small distance from the former set. So, let us fix a $t>0$ and an $\ep>0$. Apply the previous argument after~\eqref{eq:fn} to determine a $0<\delta=\delta(\ep,t)$ with $\delta<t$ and so that $\big\|\int_{t-\delta}^{t} e^{(t-s)A}\,\tilde f_n(s,v_n(s))\,ds\big\|< \ep$ for all $n\geq 1$. Then, similarly to what has been done in the last part of the proof of Proposition~\ref{prop:mild solution in K}, we write for each $n\geq 1$,
\begin{equation*}
\begin{split}
v_n(t)&= e^{tA}\,z_n+ \int_0^t e^{(t-s)A}\,\tilde f_n(s,v_n(s))\,ds \\
&= e^{\delta A}\left(e^{(t-\delta)A}\,z_n +\int_0^{t-\delta} e^{(t-s-\delta)A}\,\tilde f_n(s,v_n(s))\,ds\right) \\ &\hspace{0,4cm}+  \int_{t-\delta}^t e^{(t-s)A}\,\tilde f_n(s,v_n(s))\,ds\,.
\end{split}
\end{equation*}
Using once more Lemma~\ref{lem:bounded m bounds}, it is now easy  to determine an $R>0$ so that
\[
d\big(v_n(t),e^{\delta A}\B_R^X\big)<\varepsilon \quad\text{ for  all}\;\,n\geq 1\,,
\]
and since the operator $e^{\delta A}$ is compact and the set $\B_R^X$ is bounded, the set $e^{\delta A}\B_R^X$ is precompact and we are done.
Then, by Arzel\`{a}-Ascoli's theorem, we can assume without loss of generality that $v_n\to v$ uniformly on $[0,T]$ as $n\to\infty$.

Finally, we define $T_0$ as in \eqref{eq:T0} and we check that $v(t)$ is the mild solution for the problem for $f$ and $z$ on $[0,T_0]$. For that, it suffices to prove relation \eqref{eq:convergence}. As in the cited proof, we consider the terms named $a_n$ and $b_n$. The treatment of $b_n$ is just the same, whereas the treatment of $a_n$ is much easier now. We just need to bound
\begin{align*}
a_n&:=\left\|\, \int_0^t e^{(t-s)A}\,\tilde f_n(s,v_n(s))\,ds  - \int_0^t e^{(t-s)A}\,\tilde f(s,v_n(s))\,ds  \,\right\| \\
&\leq N \int_0^t \big\|\tilde f_n(s,v_n(s))-\tilde f(s,v_n(s))\big\|\,ds\leq N\, p_{I,j}(f_n-f) \to 0 \quad \text{as}\; n\to\infty
\end{align*}
due to the convergence of $f_n$ to $f$ in the topology $\T_B$.

The proof is finished as in Theorem~\ref{thm:CONTIN equic m bounds}, by checking that $T_0=T$.
\end{proof}
We finish the paper with a corollary that follows from the continuity results Theorems~\ref{thm:CONTIN equic m bounds} and~\ref{thm:CONTIN bounded l bounds}, combined with Theorem~\ref{thm:equivalencia topologias} on the equivalence of different topologies in some cases.
\begin{cor}
If $E\subset \LC$ has $L^1_{loc}$-bounded $l$-bounds with respect to a modulus of continuity $\widehat \theta$ in the variable $x$, then $\T_{\widetilde DD}$ is the coarsest or simplest strong topology which is good for the continuity of the mild solutions. If, in addition, $E$ has $L^1_{loc}$-equicontinuous $m$-bounds, then $\sigma_{\widetilde DD}$ is the coarsest weak topology of continuity.
\end{cor}

\end{document}